\documentclass[12pt,a4paper]{article}
\usepackage{latexsym,amssymb,amsmath,mathrsfs}
\usepackage{cite}
\usepackage{sectsty}
\usepackage[anchorcolor=blue,%
 bookmarks=true,%
 bookmarksnumbered=true,%
 dvipdfm,%
]{hyperref}
\usepackage{color}
\newtheorem{Thm}{Theorem}[section]
\newtheorem{Prop}[Thm]{Proposition}
\newtheorem{Cor}[Thm]{Corollary}
\newtheorem{Lem}[Thm]{Lemma}
\newtheorem{Defn}[Thm]{Definition}
\newtheorem{Rem}[Thm]{Remark}
\newtheorem{Ass}{Assumption}

\newenvironment{Proof}{\par\begin{trivlist}%
\item[]{\bf Proof.\ }}%
{\hfill $\square$ \end{trivlist}\par}
\newenvironment{tProof}[1]{\par\begin{trivlist}%
\item[]{\bf Proof of #1.\ }}%
{\hfill $\square$ \end{trivlist}\par}

\makeatletter \@addtoreset{equation}{section} \makeatother

\renewcommand{\P}{\mathbb{P}}
\newcommand{\E}{\mathbb{E}}
\newcommand{\R}{\mathbb{R}}
\newcommand{\N}{\mathbb{N}}

\renewcommand{\a}{\alpha}
\renewcommand{\b}{\beta}
\newcommand{\gm}{\gamma}
\newcommand{\dl}{\delta}
\newcommand{\ep}{\varepsilon}
\newcommand{\te}{\theta}
\newcommand{\lm}{\lambda}
\newcommand{\sg}{\sigma}
\newcommand{\ph}{\varphi}
\newcommand{\Dl}{\Delta}
\newcommand{\Gm}{\Gamma}
\newcommand{\Lm}{\Lambda}

\newcommand{\abs}[1]{\left| #1 \right|}
\newcommand{\abra}[1]{\left( #1 \right)}
\newcommand{\bbra}[1]{\left\{ #1 \right\}}
\newcommand{\cbra}[1]{\left[ #1 \right]}
\newcommand{\dbra}[1]{\langle #1 \rangle}
\DeclareMathOperator{\Ric}{Ric}
\DeclareMathOperator{\diam}{diam}
\DeclareMathOperator{\vol}{vol}
\DeclareMathOperator{\Cov}{Cov}
\DeclareMathOperator{\osc}{osc}
\newcommand{\wg}{\wedge}
\newcommand{\e}{\mathrm{e}}

\allsectionsfont{\normalsize\sc\center}

\title{\sf 
Monotonicity of time-dependent transportation costs 
and coupling by reflection 
}
\author{Kazumasa Kuwada\footnote{
Graduate School of Humanities and Sciences, 
Ochanomizu University, Tokyo 112-8610, Japan, 
\textit{e-mail}: \texttt{kuwada.kazumasa@ocha.ac.jp}, 
\textit{tel}: \texttt{+81-3-5978-5300}, 
\textit{fax}: \texttt{+81-3-5978-5295}
}
\ and Karl-Theodor Sturm\footnote{
Institut f\"ur Angewandte Mathematik, 
Universit\"at Bonn, 
Endenicher Allee 60, 53115 Bonn, Germany,  
\textit{e-mail}: \texttt{sturm@uni-bonn.de}
}
}
\date{}
\begin{document}
\maketitle
\raggedbottom 
\begin{abstract}
Based on a study of the coupling by reflection of 
diffusion processes, a new monotonicity in time 
of a time-dependent transportation cost between heat distribution 
is shown under Bakry-\'Emery's curvature-dimension condition
on a Riemannian manifold.
The cost function comes from the total variation 
between heat distributions on spaceforms. 
As a corollary, we obtain a comparison theorem for 
the total variation between heat distributions. 
In addition, we show that our monotonicity is stable 
under the Gromov-Hausdorff convergence of 
the underlying space 
under a uniform curvature-dimension and diameter bound. 
\end{abstract}
\begin{list}{\textbf{Keywords:}}
{\setlength{\labelwidth}{72pt}
\setlength{\leftmargin}{72pt}}
\item
transportation cost, coupling by reflection, diffusion process, 
curvature-dimension condition, total variation 
\end{list}
\begin{list}{\textbf{Mathematics Subject Classification (2010):}}
{
  \setlength{\labelwidth}{72pt}
  \setlength{\leftmargin}{72pt}
}
\item
58J65, 53C21, 60H30, 60J60, 58J35
\end{list}

\section{Introduction} 
\label{sec:intro} 

Analysis of the heat equation on manifolds or metric measure spaces 
is one of the central issues in the literature. 
Several topics such as analysis of partial differential equations, 
differential geometry and probability theory 
are interacting with each other there. 
As one of remarkable consequences of such an interaction, 
many different characterizations of 
the presence of lower Ricci curvature bound 
by means of the heat semigroup or the Brownian motion 
are revealed in \cite{Stu_Renes05}. 
Among those studies, 
recent developments in the theory of optimal transport 
enable us to interpret the heat distribution 
as a gradient curve of the relative entropy 
in the space of probability measures 
(see \cite{AGS,book_Vil2}, for instance) 
along Otto's heuristic idea in \cite{Otto_CPDE01}. 
This viewpoint provides a quite natural understanding of 
the fact that the presence of lower Ricci curvature bound 
implies a contraction property of heat distributions 
in Wasserstein distance. 
Significantly, 
this argument can bring a piece of implications 
between equivalent notions 
in \cite{Stu_Renes05} mentioned above. 
As its probabilistic counterpart, 
we can show the contraction by means of constructing 
a coupling by parallel transport of Brownian motions. 
On the other hand, there is another kind of coupling, 
called the coupling by reflection or the Kendall-Cranston coupling,  
which is also well-studied in connection with 
the Riemannian geometry of the underlying space. 
The purpose of this article is to study 
the coupling by reflection by formulating it 
in terms of the theory of optimal transport. 

To state our result, we introduce the notion of transportation cost. 
Given a function $c : M \times M \to \R$ on a state space $M$, 
a transportation cost $\mathcal{T}_c ( \mu, \nu )$ 
between two probability measures $\mu$ and $\nu$ on $M$ 
is defined as follows: 
\[
\mathcal{T}_c ( \mu , \nu ) 
: = 
\inf_{\pi \in \Pi ( \mu , \nu )} 
\int_{M \times M} c \, d \pi , 
\]
where $\Pi ( \mu , \nu )$ is the set of all couplings of 
$\mu$ and $\nu$, namely, a probability measure on $M \times M$ 
whose marginal distributions are $\mu$ and $\nu$ respectively. 
The simplest case of our result is as follows: 

\begin{Thm} \label{th:main-1}
Let $M$ be a complete Riemannian manifold with nonnegative 
Ricci curvature with the Riemannian distance $d$. 
Let us define $\ph_t (a)$ for $t , a \ge 0$ by 
\[
\ph_t (a) 
: = 
\begin{cases}
\displaystyle
\frac{1}{\sqrt{4 \pi t}} 
\int_{\R} 
\left| 
    \exp \abra{ - \frac{1}{4t} \abra{ x - \frac{a}{2} }^2 } 
    - 
    \exp \abra{ - \frac{1}{4t} \abra{ x + \frac{a}{2} }^2 } 
\right| 
dx, 
& 
t > 0, 
\\
1_{(0, \infty)} (a), 
& 
t = 0.
\end{cases}
\]
Then, for $t > 0$ and two heat distributions $\mu^{(i)}_s$ ($i=1,2$) 
generated by the Laplace-Beltrami operator, 
we have 
\begin{equation} \label{eq:mono-1}
\mathcal{T}_{\ph_{t-s_2} (d)} ( \mu^{(1)}_{s_2} , \mu^{(2)}_{s_2} )
\le 
\mathcal{T}_{\ph_{t-s_1} (d)} ( \mu^{(1)}_{s_1} , \mu^{(2)}_{s_1} )
\end{equation}
for any $0 \le s_1 \le s_2 \le t$. 
\end{Thm}

\noindent
In the full statement in Theorem~\ref{th:main2}, 
the same result as \eqref{eq:mono-1} holds 
for distributions of a diffusion process 
with an upper bound of dimension and 
a lower Ricci curvature bound 
in the sense of Bakry and \'Emery 
with an appropriate choice of $\ph_t (a)$, 
which is more complicated. 
Alternatively, \eqref{eq:mono-1} could be 
formulated as a non-expansion result of Lipschitz constants 
with respect to time-dependent metrics $\ph_t (d)$; 
see Theorem~\ref{th:gradient}. 
This allows us to interpret $\ph_t$ 
as the profile of the ``worst case'' initial data 
corresponding to $\ph_0 (d)$. 
Given $K \in \R$, 
the $L^p$-contraction in Wasserstein distance mentioned above 
means 
\begin{equation} \label{eq:W-contr}
\e^{pKt} \mathcal{T}_{d^p} ( \mu^{(1)}_t , \mu^{(2)}_t )
\le 
\e^{pKs} \mathcal{T}_{d^p} ( \mu^{(1)}_s , \mu^{(2)}_s )
\end{equation} 
for $t > s \ge 0$ and 
any heat distributions $\mu^{(i)}_s$ ($i=1,2$). 
It holds with $K = 0$ 
under the assumption in Theorem~\ref{th:main-1} 
and hence Theorem~\ref{th:main-1} 
can be regarded as an analogue of it. 
Indeed, the only difference between them 
is the choice of the cost function. 
The counterpart of Theorem~\ref{th:gradient} 
for \eqref{eq:W-contr} is the equivalence 
with Bakry-\'Emery's $L^q$-gradient estimates (see \cite{K9}). 

Let us review the history on the study of coupling by reflection, 
to explain a meaning and significance of Theorem~\ref{th:main-1}. 
We call $( X_1 (t) , X_2 (t) )$ a coupling of 
a diffusion process $X(t)$ on a state space $M$ 
if $( X_1 , X_2 )$ is a stochastic process on $M \times M$ 
and each $X_i$ behaves as $X$ on $M$ for $i =1,2$. 
The coupling by reflection on a Euclidean space, 
or the mirror coupling, 
of Brownian motions introduced in \cite{Lindvall-Rogers}
is given by the global reflection with respect to 
the hyperplane bisecting the line segment 
joining initial positions. 
With the help of Riemannian geometry, 
a coupling by reflection of 
Brownian motions on a Riemannian manifold 
is constructed by Kendall \cite{Kend} and Cranston \cite{Crans}, 
by making a coupling of their infinitesimal motions. 
In many applications, it is nice to suppose 
that they will coalesce after the coupling time, 
namely, the time when they meet. 
As a matter of fact, 
the coupling by reflection of Brownian motion 
on a Euclidean space, 
or more generally the one on a Riemannian manifold 
with nonnegative Ricci curvature, 
can meet in a finite time almost surely 
regardless of the dimension of the space. 
It is a great contrast with the case of 
observing two independent Brownian motions. 
Under a nice condition, for example, 
the presence of curvature bounds on the state space, 
this kind of coupling has provided several applications 
e.g.~in estimating the rate of convergence to equilibrium, 
functional inequalities involving heat semigroups, 
(non-)existence of harmonic maps 
(see \cite{Kend_survey} and references therein, for instance). 
As a simple example, the coupling by reflection 
under nonnegative Ricci curvature easily implies 
the Liouville property, that is, non-existence of 
nonconstant bounded harmonic functions. 
In many of those applications, we only need to know the existence of 
\emph{a coupling} $\pi$ of two distributions of $X (t)$ 
having a good estimate of $\pi ( \{ ( x, x ) \; | \; x \in M \} )$, 
which can be provided by comparing the transportation cost
in Theorem~\ref{th:main-1} at $s=t$ with that at $s=0$ 
since $\ph_0 = 1_{( 0 , \infty)}$. 
Thus the monotonicity of transportation cost in 
Theorem~\ref{th:main-1} works 
sufficiently well in applications. 

Recently, such topics as mentioned above 
is extensively studied 
on more singular metric measure spaces 
than Riemannian manifolds 
under a new, synthetic notion of curvature bounds 
(see e.g.~\cite{BBI,Lott-Vill_AnnMath09,Sturm_Ric,Sturm_Ric2}). 
However, the traditional way of studying a coupling by reflection 
of Brownian motions is based on 
the theory of stochastic differential equations 
and hence there are many difficulties to extend the original 
argument directly into analysis on such singular spaces. 
In contrast to such an approach, the statement of 
Theorem~\ref{th:main-1} completely makes sense 
even on singular spaces once we introduced 
the notion of heat distributions on it. 
Thus there seems to be some possibility 
to extend it to such cases 
though our framework in Theorem~\ref{th:main-1} or Theorem~\ref{th:main2}
is still on a Riemannian manifold. 
Actually, we obtain a partial result in this direction 
by showing that the monotonicity of the transportation cost 
stated in Theorem~\ref{th:main-1}(or Theorem~\ref{th:main2}) 
is stable under the Gromov-Hausdorff convergence of underlying spaces
(see Theorem~\ref{th:stability}) 
with a uniform curvature-dimension and diameter bound. 

One might wonder why the cost function in Theorem~\ref{th:main-1} appears. 
It is based on the fact that 
$\mathcal{T}_{\ph_{t-s} (d)} ( \mu^{(1)}_s , \mu^{(2)}_s )$ is 
a constant function in $s$ and the infimum in the definition 
of the transportation cost is attained by the coupling by reflection 
when $M$ is a Euclidean space \cite{Hsu-Stu,K4}. 
Our choice of the cost function is natural and sharp in this sense. 
Our argument in the proof of Theorem~\ref{th:main-1} 
is based on a comparison between 
the distance process for the coupling by reflection on $M$ 
and the one on a Euclidean space. 
And then the sharpness on a Euclidean space 
plays a prominent role when we deal with 
the comparison process. 
It should be remarked that, when $M$ is a Euclidean space, 
the cost function $\ph_t (d(x,y))$ coincides 
with the total variation between two heat distributions at 
time $t$ with initial distributions $\dl_x$ and $\dl_y$ respectively. 
This fact is closely related to the maximality of 
the coupling by reflection of Brownian motions on a Euclidean space 
(see \cite{Hsu-Stu,K4}, for instance). 
Indeed, by choosing the cost function as the total variation between 
heat distributions, exactly the same constancy holds on a spaceform 
since the coupling by reflection of Brownian motions is also maximal. 
As we will see, this characterization of our cost function 
leads to the following comparison theorem 
for the total variation between heat distributions. 
Let us denote the total variation norm 
by $\|\cdot \|_{\mathrm{TV}}$. 

\begin{Cor} \label{cor:TV-comparison0}
Let $M$ be a complete Riemannian manifold 
whose dimension is less than or equal to $N \in \N$ 
and whose Ricci curvature is greater than or equal to $K \in \R$.  
Then, for two heat distributions $\mu^{(1)}_s$, $\mu^{(2)}_s$ with 
$\mu^{(i)}_0 = \dl_{x_i}$ for some $x_i \in M$ ($i =1,2$), 
\begin{equation} \label{eq:TV0}
\left\| 
    \mu^{(1)}_t 
    - 
    \mu^{(2)}_t 
\right\|_{\mathrm{TV}}
\le 
\int_{\mathbb{M}_{K,N}} 
\abs{
  \tilde{p}_t^{K,N} \abra{ \tilde{x}_1 , y }
  - 
  \tilde{p}_t^{K,N} \abra{ \tilde{x}_2 , y }
}
\vol_{\mathbb{M}_{K,N}} (dy) ,
\end{equation}
where $\tilde{p}_t^{K,N} (x,y)$ is the heat kernel 
on the $N$-dimensional spaceform $\mathbb{M}_{K,N}$ 
of constant sectional curvature $K/(N-1)$ 
and $( \tilde{x}_1 , \tilde{x}_2 )$ is 
any pair of points in $\mathbb{M}_{K,N}$ 
satisfying $d( \tilde{x}_1 , \tilde{x}_2 ) = d ( x_1 , x_2 )$. 
\end{Cor}
This is a special case of Corollary~\ref{cor:TV-comparison} below. 
It seems to be natural that we can measure the total variation 
as a result of a study of the coupling by reflection 
since the coupling by reflection has been strongly related with 
estimates involving the coupling time, which yields an estimate of 
the total variation between distributions via the coupling inequality 
(see \cite{Lindv}, for instance). 

Note that, to the best of the authors' knowledge, 
an estimate of type \eqref{eq:mono-1} 
with the use of a time-dependent cost function 
is studied first in \cite[Example~4.6]{Sturm_SemigrHarm}. 
While it is discussed only on $\R^m$, 
it includes L\'evy processes as an example. 
Also note that our cost function 
in Theorem~\ref{th:main-1} (or Theorem~\ref{th:main2}) 
is a concave function of the distance function. 
It corresponds to the observation in \cite{Gang-McC_ActaMath96}, 
which says that the optimal transport map 
for a concave cost function reverses the orientation. 
Indeed, the reflection map used 
in constructing the coupling by reflection 
does so. 
Finally, we remark that 
there is a recent related result 
in \cite{Eberle_CRM11}, which studies 
a behavior of the transportation cost of 
a concave cost function 
in connection with the coupling by reflection 
of a diffusion process on $\R^m$. 

The organization of the paper is as follows. 
In the next section, 
we will give a more precise statement of our main results. 
For proving them, we will study the coupling by reflection in section 3. 
There we will follow the argument in \cite{K8,K10} in which 
we construct the coupling by reflection via an approximation 
of diffusion processes by geodesic random walks. 
It might be possible to follow an alternative approach 
in \cite{Wang_book05}. 
Section 4 is devoted to show several regularity properties 
of the function $\ph_t^{K,N} (a)$ introduced in section 2 
to describe the main theorem. 
Some explicit expressions of $\ph_t^{K,N} (a)$ 
as well as asymptotic behavior as $t \to 0$ or $t \to \infty$ 
are also given there. 
Some results in this section might be of independent interest. 
The proof of our main theorem is given in section 5. 
Though most part will follow from the result in section 3,  
we need an additional argument 
with the aid of results in section 4 to complete the proof. 
We also study new monotonicity formulae 
for time-independent transportation costs 
(Corollary~\ref{cor:mono}) 
as a consequences of the main theorem and 
results in section 4. 
In section 6, we give a short remark on gradient estimates 
for the diffusion semigroup corresponding to our main theorem. 
Though a similar gradient estimate is already 
obtained in \cite{Crans} in the same spirit, 
what we obtained is sharper in many respect. 
The stability of our main result 
under the Gromov-Hausdorff convergence 
is discussed in section 7. 
It ensures that all the results 
obtained before this section will be inherited 
to the measured Gromov-Hausdorff limit 
under a uniform curvature-dimension and diameter bound. 
In section 8, we will give a brief comment 
on the extension of results in sections 2--6 
to the time-dependent metric case. 
Note that the assumption there is satisfied 
if the metric evolves 
according to the backward Ricci flow. 

\section{Framework and the main result}
\label{sec:pre}

Let $(M, g)$ be a complete $m$-dimensional Riemannian manifold 
with $m \ge 2$. 
Let $d$ stand for the Riemannian distance on $M$. 
Let $Z$ be a $C^1$-vector field 
and we denote the generator of the form $\Dl + Z$ by $\mathcal{L}$, 
where $\Delta$ is the Laplace-Beltrami operator with respect to $g$. 
Let $( ( X(t) )_{t \in [ 0 , \infty)} , ( \P_x )_{x \in M} )$ 
be a diffusion process 
associated with $\mathcal{L}$. 
Let $(\nabla Z)^\flat$ be the symmetrization of $\nabla Z$, i.e., 
a $(0,2)$-tensor given by 
\[
(\nabla Z)^\flat (X,Y) 
  := 
\frac12 
\abra{ 
  \dbra{ \nabla_X Z , Y }
  + 
  \dbra{ \nabla_Y Z , X }
}. 
\] 
Our basic assumption is the following condition involving 
the upper dimension bound and lower Ricci curvature bound 
formulated in terms of $\mathcal{L}$: 

\begin{Ass} \label{ass:CD}
Given $K \in \R$ and $N \in [m ,\infty]$, 
the following holds: 
\[
\Ric 
- ( \nabla Z )^\flat 
- \frac{1}{N-m} Z \otimes Z 
\ge 
K g. 
\]
Here we regard the third term in the left hand side is 0 when $N = \infty$, 
and $N = m$ is permitted only when $Z \equiv 0$. 
\end{Ass}
It is well known that Assumption~\ref{ass:CD} is equivalent 
to the following curvature-dimension condition of Bakry and 
\'Emery (see e.g.~\cite{Bak97,Led_geom-Markov}): 
\[
\frac12 \abra{ 
  \mathcal{L} \dbra{ \nabla f , \nabla f } 
  - 2 \dbra{ \nabla f , \nabla \mathcal{L} f } 
}
\ge 
K \dbra{ \nabla f , \nabla f } 
+ \frac{1}{N} ( \mathcal{L} f )^2 . 
\]
This condition is equivalent to $\dim M \le N$ and $\Ric \ge K$ 
when $Z \equiv 0$. 

In order to state our main theorems, 
we introduce the notion of comparison process and 
associated transportation costs. 
Let $K \in \R$ and $N \in  [ 2, \infty ]$. 
Set $\bar{R} = \bar{R}_{K,N}$ by 

\begin{equation*}
\bar{R}_{K,N} := 
\begin{cases}
\displaystyle \sqrt{\frac{N-1}{K}} \pi & \mbox{if $K > 0$ and $N < \infty$,} 
\\ 
\infty & \mbox{otherwise.}
\end{cases}
\end{equation*}
We define $s_K$ and $c_K$ 
as a usual comparison function as follows: 
\begin{align*}
s_K (\te) 
& := 
\begin{cases} 
\displaystyle 
\frac{1}{\sqrt{K}} \sin ( \sqrt{K} \te )
& K > 0 ,
\\
\te 
& K = 0 , 
\\
\displaystyle 
\frac{1}{\sqrt{-K}} \sinh ( \sqrt{-K} \te )
& K < 0, 
\end{cases}
\\
c_K (\te) 
& := 
\begin{cases} 
\cos ( \sqrt{K} \te )
& K > 0 ,
\\
1 
& K = 0 , 
\\
\cosh ( \sqrt{-K} \te )
& K < 0 
\end{cases}
\end{align*}
and $t_K := s_K / c_K$. 
Let $\Psi = \Psi_{K,N} \: : \: ( - \bar{R}, \bar{R} )  \to \R$ 
be given by 
\begin{equation*}
\Psi_{K,N} (u)
: = 
\begin{cases}
\displaystyle 
- 2 K t_{K/(N-1)} \abra{ \frac{u}{2} } 
& \mbox{if $N < \infty$,} 
\\
\displaystyle 
- K u 
& \mbox{otherwise.}
\end{cases}
\end{equation*}
Let us define a diffusion process 
$\rho (t) = \rho_{K,N,a} (t)$, $t \ge 0$
on $\overline{( - \bar{R} , \bar{R} )} \subset \R$ 
as a solution to the following stochastic differential equation: 
\begin{align} \label{eq:rho0}
d \rho_{K,N,a} (t) 
& = 
2 \sqrt{2} d \b (t) 
+ \Psi ( \rho_{K,N,a} (t) ) d t , 
\\ \nonumber
\rho_{K,N,a} (0) 
& = 
a. 
\end{align}
Note that, when $\bar{R} < \infty$, 
both $-\bar{R}$ and $\bar{R}$ are entrance boundary for $\rho(t)$. 
For $t \ge 0$, let us define 
$\ph^{K,N}_t \, : \, \overline{[ 0 , \bar{R} )} \to [0,1]$ 
by 
\begin{equation*}
\ph^{K,N}_t ( a ) 
: = 
\P \cbra{ \inf_{0 \le s \le t} \rho_{K,N,a} (s) > 0 }. 
\end{equation*}

\begin{Rem} \label{rem:comparison}
\begin{enumerate}
\item \label{comparison}
The process $\rho_{K,N,a}$ comes from 
the coupling by reflection on the spaceform. 
Actually, when $N \in \N$, 
a simple computation implies that 
the distance process $d ( \mathbf{X} (t) )$ 
for the coupling by reflection 
$\mathbf{X} (t) = ( X_1 (t) , X_2 (t) )$
of Brownian motions on the spaceform $\mathbb{M}_{N,K}$ 
solves the stochastic differential equation defining 
$\rho_{K,N,a}$ with $a = d ( \mathbf{X} (0) )$. 

\item \label{reflection}
Since $-\rho_{K,N,a}$ has the same law as $\rho_{K,N,-a}$, 
the reflection map $x \mapsto -x$ on 
$\overline{(- \bar{R} , \bar{R})}$ provides 
a so-called `reflection structure' in \cite{K4}. 
It is shown in \cite{K4} that 
the mirror coupling 
for $\rho_{K,N,a}$ and $\rho_{K,N,-a}$ is maximal 
in such a case. 
As a result, we have 
\[
\ph^{K,N}_t (a) 
= 
\left\| 
    \P \circ ( \rho_{K,N,a} (t) / 2 )^{-1} 
    - 
    \P \circ ( \rho_{K,N,-a} (t) / 2 )^{-1} 
\right\|_{\mathrm{TV}}. 
\]
In particular, we can easily verify that 
$\ph^{0,N}_t$ equals to $\ph_t$ in Theorem~\ref{th:main-1}. 
Moreover, 
$\E [ \ph_{t-s} ( | \rho (s) | ) ]$ 
is a constant function in $s \in [ 0, t ]$ 
by \cite[Lemma~3.4]{K4}. 

\item \label{spaceform}
When $N \in \N$, the coupling $\mathbf{X} (t)$ 
by reflection of 
Brownian motions $B (t)$ on $\mathbb{M}_{K,N}$ 
is maximal by the same reasoning 
(see \cite[Theorem~5.1 and Example~4.6]{K4}). 
Thus we have 
\[
\ph^{K,N}_t (a) 
= 
\left\| 
    \P_{\tilde{x}_1} \circ B (t)^{-1} 
    - \P_{\tilde{x_2}} \circ B (t)^{-1} 
\right\|_{\mathrm{TV}} 
\]
for any pair of points $( \tilde{x}_1 , \tilde{x}_2 )$ 
in $\mathbb{M}_{K,N}$ 
satisfying $d ( \tilde{x}_1 , \tilde{x}_2 ) = a$ 
and 
$\E [ \ph_{t-s} ( d ( \mathbf{X} (s) ) ) ]$ 
is a constant function in $s \in [ 0, t ]$. 
In particular, the right hand side of 
\eqref{eq:TV0} equals to $\ph^{K,N}_t ( d (x_1 , x_2 ))$. 
\end{enumerate}
\end{Rem}

Now we are in turn 
to state our first main theorem as follows: 
\begin{Thm} \label{th:main}
Suppose that Assumption~\ref{ass:CD} holds. 
Then, for any $x_1 , x_2 \in M$, 
there exists a coupling 
$\mathbf{X}(t) = ( X_1 (t) , X_2 (t) )_{t \ge 0}$ of 
$\mathcal{L}$-diffusion processes 
starting from $( x_1 , x_2 )$ 
such that, for any $t > 0$ and $s \ge 0$, 
\begin{equation*}
\E \cbra{ 
  \ph^{K,N}_t ( d ( \mathbf{X} (s) ) )
}
\le 
\ph^{K,N}_{t+s} d ( x_1 , x_2 ). 
\end{equation*}
\end{Thm}
Indeed, as we will see, 
a coupling $\mathbf{X} (t)$ appeared in Theorem~\ref{th:main} 
will be given as the coupling by reflection. 
Theorem~\ref{th:main} yields the following corresponding property 
described in terms of $\mathcal{T}_{\ph_t ( d )}$. 
This is our second main theorem: 
\begin{Thm} \label{th:main2}
Suppose that Assumption~\ref{ass:CD} holds. 
For $i=1, 2$ and $\mu^{(i)} \in \mathcal{P} (M)$, 
let $\mu^{(i)}_t$ be the distribution of $X(t)$ 
with the initial distribution $\mu^{(i)}$. 
Then, for any $t > 0$, 
$\mathcal{T}_{\ph_{t-s} ( d )} ( \mu^{(1)}_s , \mu^{(2)}_s )$ is 
a nonincreasing function of $s \in [ 0 , t ]$. 
That is, for $0 \le s_1 \le s_2 \le t$, 
\begin{equation} \label{eq:mono1}
\mathcal{T}_{\ph_{t-s_2} ( d )} ( \mu^{(1)}_{s_2} , \mu^{(2)}_{s_2} ) 
\le
\mathcal{T}_{\ph_{t-s_1} ( d )} ( \mu^{(1)}_{s_1} , \mu^{(2)}_{s_1} ) . 
\end{equation} 
\end{Thm}
As a result of Theorem~\ref{th:main2}, 
we can compare 
$\mathcal{T}_{\ph_{t - s} ( d )} ( \mu^{(1)}_s , \mu^{(2)}_s )$
at $s = t$ with the one at $s = 0$ 
to obtain an estimate of the total variation 
between distributions of the diffusion process $X (t)$. 
In particular, when $\mu^{(1)}_0$ and $\mu^{(2)}_0$ are Dirac measures, 
we obtain the following comparison theorem 
thanks to Remark~\ref{rem:comparison} \ref{reflection}: 
\begin{Cor} \label{cor:TV-comparison}
Suppose that Assumption~\ref{ass:CD} holds. 
Then, for $x_1 , x_2 \in M$ and $t > 0$, 
\begin{multline*}
\left\| 
    \P_{x_1} \circ X (t)^{-1} 
    - 
    \P_{x_2} \circ X (t)^{-1} 
\right\|_{\mathrm{TV}}
\\
\le 
\left\| 
    \P \circ ( \rho_{K,N,d (x_1 , x_2)} (t) /2 )^{-1} 
    - 
    \P \circ ( \rho_{K,N, -d (x_1 , x_2)} (t) /2 )^{-1} 
\right\|_{\mathrm{TV}}. 
\end{multline*}
\end{Cor}
When $N \in \N$, it immediately implies Corollary~\ref{cor:TV-comparison0} 
by virtue of Remark~\ref{rem:comparison} \ref{spaceform}. 

Note that, by taking $t \to \infty$ in \eqref{eq:mono1} 
after a suitable rescaling, 
we can obtain a similar monotonicity formula 
whose cost is independent of $t$. 
See Corollary~\ref{cor:mono} below. 
Especially, when $K < 0$, 
it does not seem to be known in the literature. 
\begin{Rem} \label{rem:BMC} 
When $K > 0$, it is shown in \cite{K11} that 
under Assumption~\ref{ass:CD} 
the Bonnet-Myers type diameter bound 
\[
\diam (M) \le \pi \sqrt{\frac{N-1}{K}}
\] 
holds. 
Moreover, the equality holds only 
when $N = m$, $Z \equiv 0$ and $M$ is isometric to 
$N$-dimensional sphere of 
constant sectional curvature $K/(N-1)$. 
In the case of equality, 
the assertion in Theorem~\ref{th:main} is obvious 
by Remark~\ref{rem:comparison} \ref{spaceform} 
and hence we may assume $\diam (M) < \pi \sqrt{(N-1)/K}$ 
in the sequel. 
\end{Rem}

\section{Proof of \protect{Theorem~\ref{th:main}}}
\label{sec:pf}

We will show that the coupling by reflection 
studied in \cite{K8} (cf.~\cite{K10}) 
satisfies the assertion of Theorem~\ref{th:main}
under Assumption~\ref{ass:CD}. 
We begin with reviewing the construction of 
the coupling by reflection. 
Let $( \xi_n )_{n \in \N}$ be 
independent random variables 
all of which are uniformly distributed 
on the unit disk on $\R^m$. 
Let $( \gm_{xy} )_{x,y \in M}$ be 
a measurable family of 
unit-speed minimal geodesics 
defined on $[0, d (x,y)]$ 
such that $\gm_{xy}$ joins $x$ and $y$. 
Without loss of generality, 
we may assume that 
$\gm_{xy}$ is symmetric, that is, 
$
\gm_{xy} ( d(x,y) - s ) 
= 
\gm_{yx} (s)
$ 
holds. 
Let us define 
$\tilde{m}_{xy} \: : \: T_x M \to T_x M$ 
by 
\[
\tilde{m}_{xy} v 
:= 
v 
- 2 \dbra{ v, \dot{\gm}_{xy} (0) } 
\dot{\gm}_{xy} (0)
. 
\] 
This is a reflection with respect to a hyperplane 
which is perpendicular to $\dot{\gm}_{xy}$. 
Let $/\!\!/_{\gm}$ be the parallel transport 
along a curve $\gm$. 
Let us define $m_{xy} \: : \: T_x M \to T_y M$ by 
$m_{xy} := /\!\!/_{\gm_{xy}} \circ \tilde{m}_{xy}$. 
Clearly $m_{xy}$ is an isometry. 
Set $D (M) : = \{ (x,x) \; | \; x \in M \}$.  
Let $\Phi \: : \: M \to \mathscr{O} (M)$ 
be a measurable section of 
the orthonormal frame bundle 
$\mathscr{O} (M)$ of $M$. 
Let us define two measurable maps 
$\Phi_i \: : \: M \times M \to \mathscr{O} (M)$ 
for $i=1,2$ by 
\begin{align*}
\Phi_1 (x,y) & 
:= \Phi (x), 
\\ 
\Phi_2 (x,y) 
& := 
\begin{cases}
m_{xy} \Phi_1 (x,y), 
& 
(x,y) \in M \times M \setminus D(M), 
\\ 
\Phi (x), 
& 
(x,y) \in D (M).
\end{cases}
\end{align*}
Take $x_1 , x_2 \in M$. 
Let $t_n^{\a} : = \a^2 n$ for $n \in \N_0$. 
By using $\Phi_i$, 
we define a coupled geodesic random walk 
$\mathbf{X}^\a (t) = ( X_1^\a (t) , X_2^\a (t) )$ 
with a scale parameter $\a$ 
by $X^\a_i (0) = x_i$ and, 
for $t \in [ t_n^{\a} , t_{n+1}^{\a} ]$, 
\begin{align} 
\nonumber
\tilde{\xi}_{n+1}^i
& : = 
\sqrt{2(m+2)} 
\Phi_i
\abra{ 
  \mathbf{X}^\a ( t_n^{\a} ) 
} 
\xi_{ n + 1 } ,
\\ \nonumber 
X_i^\a ( t ) 
& := 
\exp_{X_i^\a ( t_n^{\a} )}
\bigg(
  \frac{ t - t_n^{\a} }{\a^2}  
  \Big( 
    \a \tilde{\xi}_{n+1}^i
     + 
    \a^2 Z 
  \Big)
\bigg)  
\end{align}
for $i = 1, 2$, 
where $\exp_x$ is 
the exponential map at $x$. 
Let us denote $C ( [ 0, \infty ) \to M \times M)$ 
and 
$C ( [ 0 , \infty ) \to [ - \bar{R} , \bar{R} ] )$ 
equipped with the topology of 
compact uniform convergence 
by $\mathscr{C}$ and $\mathscr{C}_1$ respectively. 

In what follows, we assume Assumption~\ref{ass:CD}. 
Then, by \cite[Theorem~3.1]{K10} (also see references therein), 
$X_i^\a (t)$ converges in law in $C ( [ 0 , \infty ) \to M )$ 
to an $\mathcal{L}$-diffusion process 
starting from $x_i$ for $i = 1,2$ respectively. 
Thus $( \mathbf{X}^\a )_{\a > 0}$ is tight and 
hence a subsequential limit 
$\mathbf{X}^{\a_k} \to 
\mathbf{X} = ( X_1 , X_2 )$ 
in law in $\mathscr{C}$ exists. 
We fix such a subsequence $( \a_k )_{k \in \N}$. 
In the rest of this paper, 
we use the same symbol $\mathbf{X}^\a$ 
for the subsequence $\mathbf{X}^{\a_k}$ 
and the term ``$\a \to 0$'' always means 
the subsequential limit ``$\a_k \to 0$''. 
Let $\tau^*$ be the first hitting time 
to $D(M)$ of $\mathbf{X}$. 
Then we define \emph{a coupling by reflection} 
$\mathbf{X}^* = ( X_1^* , X_2^* )$ by 
\begin{equation*}
\mathbf{X}^* (t) 
: = 
\begin{cases} 
\mathbf{X} (t) 
& \mbox{if $t < \tau^*$,}
\\
( X_1 (t) , X_1 (t) ) 
& 
\mbox{if $t \ge \tau^*$.}
\end{cases}
\end{equation*}
Since $\tau^*$ is a stopping time with respect to 
the filtration generated by $\mathbf{X}$, and 
$X_i$ ($i=1,2$) is a solution to the martingale problem 
associated with the same filtration, $\mathbf{X}^*$ 
is again a coupling of $\mathcal{L}$-diffusion process. 

Fix a reference point $o \in M$. 
For $R > 0$, let 
$
\sg_R 
  \: : \: 
\mathscr{C}_1 \to 
[ 0 , \infty ]
$ 
be given by 
$
\sg_R (w) 
  : = 
\inf \bbra{ 
  t \in [ 0 , \infty ) 
  \; | \;
  w(t) \ge R 
}
$. 
We define $\hat{\sg}_R^i$ ($i=1,2$) and 
$\hat{\sg}_R$ by 
$
 \hat{\sg}_R^{i} 
 : = 
 \sg_R ( d ( o , X_i^\a (\cdot) ) )
$
and 
$\hat{\sg}_R : = \hat{\sg}_R^1 \wg \hat{\sg}_R^2$. 
Proposition~3.4 in \cite{K10} says that 
\begin{equation} \label{eq:cons}
\lim_{R \to \infty} 
\limsup_{\a \to 0} 
\P [ \hat{\sg}_R  < \infty ] = 0 
\end{equation}
holds. 

We next review a difference inequality of 
$d ( \mathbf{X}^\a (t) )$. 
To describe it, 
we will introduce some notations. 
For simplicity of notations, 
let us denote 
$\gm_{X_1^\a ( t_n^{\a} ) X_2^\a ( t_n^{\a} )}$, 
$m_{X_1^\a ( t_n^{\a} ) X_2^\a ( t_n^{\a} )}$ 
and $d ( \mathbf{X}^\a ( t_n^\a ) )$  
by $\gm_n$, $m_n$ and $r^\a (n)$ respectively. 
Let $\tilde{\xi}_{n+1}^\perp (0)$ be the orthogonal projection 
of $\tilde{\xi}_{n+1}^1$ to the hyperplane being perpendicular 
to $\dot{\gm}_n (0)$, 
that is, 
$2 \tilde{\xi}_{n+1}^\perp (0) := ( 1 + m_n ) \tilde{\xi}_{n+1}^1$. 
We denote a vector field along $\gm_n$ given 
by parallel transport of $\xi_{n+1}^\perp (0)$ by 
$( \xi_{n+1}^\perp (s) )_{s \in [0, r^\a (n)]}$. 
Let us define a weight function $h_{n+1} = h^{K,N}_{n+1}$ 
on $[0, r^\a (n) ]$ and 
a vector field $V_{n+1} = V_{n+1}^{K,N}$ 
along $\gm_n$ by 
\begin{align*}
h^{K,N}_{n+1} (s) 
& := 
\begin{cases}
\displaystyle 
c_{K/(N-1)} \abra{ 
  \frac{r^\a (n) }{2} 
}^{-1}
c_{K/(N-1)} \abra{ 
  \abra{ \frac{s - r^\a (n)}{2} } 
}
& \mbox{if $N < \infty$,}
\\
1 & \mbox{if $N = \infty$,} 
\end{cases}
\\
V_{n+1}^{K,N} (s)
& : = 
h_{n+1} (s) \tilde{\xi}_{n+1}^\perp (s) . 
\end{align*}
Recall that 
we are assuming $\diam (M) < \pi \sqrt{(N-1)/K}$ 
when $K > 0$ and $N < \infty$ 
(see Remark~\ref{rem:BMC}).  
Hence $h_{n+1}$ is well-defined. 
For a smooth curve $\gm$ and 
vector fields $V$ and $W$ along $\gm$, 
we denote the index form by $I_\gm (V, W)$. 
When $V=W$, 
we use the symbol $I_\gm (V)$ 
for $I_{\gm} (V,W)$. 
Take $v \in \R^m$. 
Let us define $\lm_{n+1}$ and $\Lm_{n+1}$ 
by 
\begin{align*}
\lm_{n+1} 
& : = 
\begin{cases}
2 \sqrt{2}
\dbra{ 
  \tilde{\xi}_{n+1}^1 (0) , 
  \dot{\gm}_n (0) 
}
& \mbox{if $\mathbf{X}^\a (t_n^\a) \notin D(M)$}, 
\\
2\sqrt{2} \sqrt{m+2} 
\dbra{ \xi_{n+1} , v } 
& \mbox{otherwise}, 
\end{cases}
\\
\Lm_{n+1} 
& : = 
\Bigg( 
\left. 
  \dbra{ Z ( t_n^{\a} ) , \dot{\gm}_n (s) }
\right|_{s=0}^{r^\a (n)}
+ 
I_{\gm_n} \abra{ V_{n+1} } 
\Bigg) 1_{\{ \mathbf{X}^\a ( t_n^{\a} ) \notin D(M) \} }. 
\end{align*}
For $\dl \ge 0$, let us define 
$\tau_\dl \: : \: \mathscr{C}_1 \to [ 0 , \infty ]$ 
by $\tau_\dl (w) := \inf \bbra{ t \ge 0 \; | \; w(t) \le \dl }$. 
We also define  
$
\hat{\tau}_\dl 
$
by 
$
\hat{\tau}_\dl 
: = 
\tau_\dl ( d ( \mathbf{X}^a (\cdot) ) )
$. 
In the sequel, we fix $\dl \in (0,1)$ and $R > 1$. 
The first goal is to prove the following difference inequality 
for $r^\a (n)$: 
\begin{Prop} \label{prop:DI}
For each $\ep > 0$, there exists a family of events 
$E_\ep^\a$ with $\lim_{\a \to 0} \P [ E_\ep^\a ] = 1$
such that 
\begin{equation*}
r^\a (n+1) 
  \le 
r^\a (n) 
+ \a \lm_{n+1} 
+ \a^2 \Psi ( r^\a (n) ) 
+ \ep \a^2 
\end{equation*}
holds for $n \in \N$ 
with $t_n^\a < \hat{\tau}_\dl \wg \hat{\sg}_{R}$ 
on $( E_\ep^\a )^c$ for sufficiently small $\a$. 
\end{Prop}
We will prove this assertion 
by a similar argument as in \cite{K8,K10}. 
Thus we only give a brief sketch of arguments. 
It consists of the following three lemmata. 
The following is shown in the same way 
as \cite[Lemma~3]{K8} or \cite[Lemma~4.4]{K10} 
by using the second variation formula of arclength 
with a careful treatment of singularities 
arising from the cutlocus. 
\begin{Lem} \label{lem:c-2var}
For $n \in \N_0$, 
we have 
\begin{equation} \label{eq:c-2var}
r^\a (n+1) 
 \le 
r^\a (n) + \a \lm_{n+1} + \a^2 \Lm_{n+1} + o (\a^2 ) 
\end{equation}
when 
$n < \hat{\tau}_\dl \wg \hat{\sg}_{R}$ 
and $\a$ is sufficiently small. 
Moreover, we can control the error term $o(\a^2)$ 
uniformly in the position of $\mathbf{X}^\a$. 
\end{Lem}
Set $\mathscr{F}_n : = \sg ( \xi_1 , \ldots , \xi_n )$ and 
$\bar{\Lm}_{n+1} : = \E \cbra{ \left. \Lm_{n+1} \right| \mathscr{F}_n }$. 
For $\ep > 0$ and $R > 0$, 
let us define an event $\tilde{E}_\ep^\a$ by 
\begin{equation*}
\tilde{E}_\ep^\a 
: = 
\bbra{
  \sup_{
        t_n^\a \le \hat{\sg}_{R}
  }
  \sum_{j=1}^n 
  ( \Lm_j - \bar{\Lm}_j ) 
  \le \frac{\ep}{2 \a^2}
}. 
\end{equation*}
By following arguments 
in \cite[Lemma~6]{K8} or \cite[Lemma~4.5]{K10} 
which are based on the Doob submartingale inequality, 
we obtain the following. 
\begin{Lem} \label{lem:LLN} 
For any $\ep > 0$ and $R > 0$, 
$\P [ \tilde{E}_\ep^\a ]$ tends to $1$ as $\a \to 0$. 
\end{Lem}
Lemma~\ref{lem:LLN} ensures to replace $\Lm_{n+1}$ in Lemma~\ref{lem:c-2var} 
with $\bar{\Lm}_{n+1}$ with small errors on $( \tilde{E}_\ep^\a )^c$. 
Thus, the proof of Proposition~\ref{prop:DI} will be completed 
with $E_\ep^\a = \tilde{E}_\ep^\a$ 
once we show the following: 
\begin{Lem} \label{lem:Exp}
$\bar{\Lm}_{n+1} \le \Psi ( r^\a (n) )$. 
\end{Lem}

\begin{Proof}
Note that we have 
\begin{align} \nonumber
\left.  
  \dbra{ Z , \dot{\gm}_n (s) }
\right|_{s=0}^{r^\a (n)}
& = 
\left.
  h_{n+1} (s)^2 
  \dbra{ Z , \dot{\gm}_n (s) }
\right|_{s=0}^{r^\a (n)}
\\ \nonumber 
& = 
\int_0^{r^\a (n)} 
\Big(
  h_{n+1} (s)^2 ( \nabla Z )^\flat 
  ( \dot{\gm}_n (s) , \dot{\gm}_n (s) ) 
\\ \label{eq:Exp1}
& \hspace{12em}
    + 
  2 h_{n+1}' (s) h_{n+1} (s) 
  \dbra{ Z , \dot{\gm}_n (s) }
\Big)
d s . 
\end{align}
By an easy computation, we obtain 
$\E [ \xi_1 ] = 0$ and 
$\Cov (\sqrt{2(m+2)}\xi_1) = 2 \mathrm{Id}$. 
Thus we have
\begin{equation} \label{eq:Exp2}
I_{\gm_n}
\abra{ V_{n+1} } 
 = 
\int_0^{r^\a (n)} 
  \abra{
    (m-1) h_{n+1}' (s)^2 
    - 
    \Ric ( \dot{\gm}_n (s) , \dot{\gm}_n (s) ) 
    h_{n+1} (s)^2 
  }
d s . 
\end{equation}
Combining \eqref{eq:Exp1} and \eqref{eq:Exp2} 
with the definition of $\bar{\Lm}_n$, 
we obtain 
\begin{align} \nonumber
\bar{\Lm}_{n+1} 
& : = 
\Bigg( 
\int_0^{r^\a (n)} 
\Big(
    h_{n+1} (s)^2 ( \nabla Z )^\flat 
  \abra{ 
    \dot{\gm}_n (s) , 
    \dot{\gm}_n (s) 
  }
    + 
  2 h_{n+1}' (s) h_{n+1} (s) 
  \dbra{ Z , \dot{\gm}_n (s) }
\\ \label{eq:Exp3}
& \hspace{3em} + 
  (m-1) h_{n+1}' (s)^2 
    - 
  \Ric
  \abra{ 
    \dot{\gm}_n (s) , 
    \dot{\gm}_n (s) 
  } 
  h_{n+1} (s)^2 
\Big)
d s  
\Bigg) 1_{\{ \mathbf{X}^\a ( t_n^{\a} ) \notin D(M) \} }. 
\end{align}
Thus, when $N = \infty$, 
the conclusion easily follows 
from Assumption~\ref{ass:CD}. 
When $N = m$, 
$Z \equiv 0$ holds and 
Assumption~\ref{ass:CD} 
means $\Ric \ge K$. 
Thus an easy computation in \eqref{eq:Exp3} 
yields the conclusion. 
When $m < N < \infty$, 
the arithmetic geometric mean inequality 
implies 
\begin{multline}
2 h_{n+1}' (s) h_{n+1} (s) 
  \dbra{ Z , \dot{\gm}_n (s) }
 \le 
(N-m) h_{n+1}' (s)^2 
 + 
\frac{1}{N-m} h_{n+1} (s)^2 
\dbra{ Z , \dot{\gm}_n (s) }^2 
\\ \label{eq:Exp4}
 = 
(N-m) h_{n+1}' (s)^2 
 + 
\frac{1}{N-m} h_{n+1} (s)^2 
Z \otimes Z 
( \dot{\gm}_n (s) , \dot{\gm}_n (s) ). 
\end{multline} 
By substituting \eqref{eq:Exp4} into \eqref{eq:Exp3}, 
we obtain 
\begin{align*}
\bar{\Lm}_{n+1} 
& \le 
\Bigg( 
\int_0^{r^\a (n)} 
\Big(
  (N-1) h_{n+1}' (s)^2 
\\ 
& \quad
  + h_{n+1} (s)^2 
  \abra{
    \frac{1}{N-m} Z \otimes Z 
      +
    ( \nabla Z )^\flat 
      - 
    \Ric 
  } 
  \abra{ 
    \dot{\gm}_n (s) , 
    \dot{\gm}_n (s) 
  }
\Big)
d s  
\Bigg) 1_{\{ \mathbf{X}^\a ( t_n^{\a} ) \notin D(M) \} }. 
\end{align*}
Hence Assumption~\ref{ass:CD} 
reduces the assertion 
to the same computation as in the case $N = m$. 
\end{Proof}

Set $a := d (x_1 , x_2)$. 
Let $\rho_{K,N,a}^\a (t)$ be a discrete approximation of 
$\rho_{K,N,a} (t)$ defined inductively 
by $\rho_{K,N,a}^\a (0) = a$ and 
for $t \in [ t_n^\a , t_{n+1}^\a ]$ 
\begin{equation*}
\rho_{K,N,a}^\a (t) 
: = 
\rho_{K,N,a}^\a (t_n^\a) 
  + 
\frac{t - t_n^\a}{\a^2} 
\abra{ 
  \a \lm_{n+1} 
    + 
  \a^2 \Psi ( \rho_{K,N,a}^\a (t_n^\a) )
}. 
\end{equation*} 
For $R > 0$, let us define $R^* > 1$ by 
\begin{equation*}
R^* := 
\begin{cases}
\displaystyle 
\bar{R} - \frac{1}{R} 
& \mbox{if $K > 0$ and $N < \infty$,}
\\ 
\displaystyle 
R 
& \mbox{otherwise.} 
\end{cases}
\end{equation*}
The following comparison theorem is crucial 
for the proof of Theorem~\ref{th:main}. 

\begin{Prop} \label{prop:compare}
For $T > 0$, $R > 0$ and $\ep > 0$, 
there exists a constant $C ( \ep , T ) \ge 0$ 
satisfying $\lim_{\ep \to 0} C ( \ep , T ) = 0$ 
such that 
\begin{equation*}
d ( \mathbf{X}_t^\a ) 
  \le 
\rho_{K,N,a}^\a (t) + C ( \ep , T )
\end{equation*}
holds for 
$
t < \hat{\tau}_\dl 
\wg \hat{\sg}_{R} 
\wg \sg_{R^*} ( \rho_{K,N,a}^\a ) 
\wg T
$ 
on $( E_\ep^\a )^c$ for sufficiently small $\a$. 
\end{Prop}

\begin{Proof}
By \cite[Corollary~3.6(i)]{K10}, 
it suffices to show the assertion 
only when $t = t_n^\a$ for some $n \le n^{(\a)}$ 
(cf. \cite[Lemma~3.10]{K10}). 
For simplicity of notations, 
we denote $\rho_{K,N,a}^\a ( t_n^\a )$ 
by $\rho^\a (n)$. 
Applying Proposition~\ref{prop:DI}, we obtain 
\begin{equation} \label{eq:compare1}
r^\a (n+1) - \rho^\a ( n+1 ) 
  \le 
r^\a (n) - \rho^\a ( n )
  + 
\a^2 ( \Psi ( r^\a (n) ) - \Psi ( \rho^\a (n) ) )
+ \ep \a^2 
\end{equation}
for $n \le n^{(\a)}$ 
with 
$
t_n^\a < 
\hat{\tau}_\dl \wg \hat{\sg}_{R} \wg \sg_{R^*} (\rho_{K,N,a}^\a)
$ 
on $E_\ep^\a$. 
Under our assumption on $t = t_n^\a$, 
$r^\a (n) \in [ \dl, R ]$ and 
$\rho^\a (n) \in [ 0 , R^* ]$ hold. 
Note that $\Psi$ is bounded 
on $[ 0 , \diam (M) \wg R^* ]$. 
Let $f_\a \: : \: \R \to \R$ 
be a function of class $C^2$ satisfying 
the following conditions: 
\begin{enumerate}
\item \label{positive}
$f_\a (x) = 0$ for $x \le 0$ 
and 
$f_\a (x) = x + \a / 2$ for $x \ge \a$,
\item \label{convex}
$f_\a$ is convex, 
\item \label{uniform}
$\displaystyle \limsup_{\a \to 0} \a^2 \sup_{u \in \R} f''_\a (u) < C$ 
for some $C > 0$ 
\end{enumerate}
(cf.~the proof of \cite[Lemma~3.10]{K10}). 
By \eqref{eq:compare1}, 
the Taylor expansion 
together with the condition \ref{uniform} of $f_\a$ 
yields 
\begin{multline} \label{eq:compare2}
f_\a ( r^\a (n+1) - \rho^\a (n+1) ) 
  \le 
f_\a ( r^\a (n) - \rho^\a (n) )
\\
  + 
\a^2 f_\a' ( r^\a (n) - \rho^\a (n) )
\abra{ \Psi (r^\a (n)) - \Psi (\rho^\a (n)) } 
+ 2 \ep \a^2  
\end{multline}
for sufficiently smaller $\a$ than $\ep$. 
Since $\Psi$ is nonincreasing, 
properties \ref{positive} and \ref{convex}
of $f_\a$ imply 
\[
f_\a' ( r^\a (n) - \rho^\a (n) )
\abra{ \Psi (r^\a (n)) - \Psi (\rho^\a (n)) } 
\le 0 . 
\]
Thus, an iteration of \eqref{eq:compare2} 
together with the fact $f_\a (x) + \a / 2 \ge x \vee 0$ 
yield 
\begin{equation*}
( r^\a (n) - \rho^\a (n) )_+  
  \le 
f_\a ( r^\a (n) - \rho^\a (n) ) 
+ \frac{\a}{2} 
  \le 
2 \ep \a^2 n  + \ep  
\end{equation*}
for $\a \le 2 \ep$. 
Since $t_n^\a = \a^2 n \le T$, 
the conclusion follows. 
\end{Proof}

Now we are in position to give a crucial step of the proof of 
Theorem~\ref{th:main}. 

\begin{Prop} \label{prop:main0}
For any nondecreasing continuous function 
$\psi \: : \: \overline{[  0 , \bar{R} )} \to [ 0, 1 ]$ 
with $\psi (0) = 0$, 
we have
\[
\E \cbra{
  \psi ( d ( \mathbf{X}^* (s) ) )
}
\le 
\E \cbra{ 
  \psi ( \rho (s) ) 
  \; ; \; 
  \tau_0 ( \rho ) > s 
} . 
\]
\end{Prop}
\begin{Proof}
Take $\dl > 0$, $R > 1$ and $t > s$. 
Let $\ep > 0$ be so small that 
$C ( \ep , t ) < \dl / 2$. 
By virtue of Proposition~\ref{prop:DI}, 
for sufficiently small $\a$, 
\begin{multline} \label{eq:loc1}
\E [ \psi ( d ( \mathbf{X}^\a (s) ) ) ] 
\le 
\E \cbra{ 
  \psi ( d ( \mathbf{X}^\a (s) ) )
  \; ; \; 
  \{ \hat{\tau}_\dl > s \} 
    \cap
  \{ \hat{\sg}_R > s \} 
    \cap 
  ( E_\ep^\a )^c
}
\\
  + 
\P [ \hat{\sg}_R \le s ]
  + 
\E \cbra{ 
  \psi ( d ( \mathbf{X}^\a (s) ) )
  \; ; \; 
  \hat{\tau}_\dl \le s 
} 
  + 
\ep .
\end{multline}
By Proposition~\ref{prop:compare} and the choice of $\ep$, 
\begin{multline} \label{eq:rep1}
\E \cbra{ 
  \psi ( d ( \mathbf{X}^\a (s) ) )
  \; ; \; 
  \{ \hat{\tau}_\dl > s \} 
    \cap
  \{ \hat{\sg}_R > s \} 
    \cap 
  ( E_\ep^\a )^c
}
\\
\le 
\E \cbra{ 
  \psi ( \rho^\a (s) + C ( \ep , t ) ) 
  \; ; \; 
  \tau_{\dl / 2} ( \rho^\a ) \wg \sg_{R^*} ( \rho^\a ) > s 
}  
+ 
\P [ \sg_{R^*} ( \rho^\a ) \le s ] . 
\end{multline}
Let us define $\tilde{\Psi} \: : \: [ 0, \infty ) \to \R$ by 
\[
\tilde{\Psi} (u) 
  := 
( \Psi (u) \wg | \Psi ( (2R)^* ) | ) \vee ( - | \Psi ( - (2R)^* ) | )
. 
\]
We define $\tilde{\rho}^\a$ and $\tilde{\rho}$ 
by replacing $\Psi$ with $\tilde{\Psi}$ 
in the definition of $\rho^\a$ and $\rho$ respectively. 
Since $\tilde{\Psi}(u) = \Psi (u)$ for $u \in [ 0 , R^* ]$, 
we obtain 
\begin{align} \nonumber 
\E \big[ 
  \psi ( \rho^\a (s) + C ( \ep , t ) ) 
&
  \; ; \; 
  \tau_{\dl / 2} ( \rho^\a ) \wg \sg_{R^*} ( \rho^\a ) > s 
\big]  
\\ \label{eq:rep2}
& \le  
\E \cbra{ 
  \psi ( \tilde{\rho}^\a (s) + C ( \ep , t ) ) 
  \; ; \; 
  \tau_{\dl / 2} ( \tilde{\rho}^\a ) > s 
},   
\\ \label{eq:rep3}
\P [ \sg_{R^*} ( \rho^\a ) \le s ] 
& = 
\P [ \sg_{R^*} ( \tilde{\rho}^\a ) \le s ] . 
\end{align}
Since $\tilde{\Psi}$ is bounded and continuous, 
we can easily show that $\tilde{\rho}^\a$ converges 
in law to $\tilde{\rho}$ in $C ( [ 0 , \infty ) \to \R )$. 
Note that the following holds:  
\begin{equation*}
\overline{
  \bbra{ 
    w \in \mathscr{C}
    \; ; \; 
    \tau_{\dl / 2 } (d (w) )  > s 
  }
}
\subset 
\bbra{ 
  w \in \mathscr{C}
  \; ; \; 
  \tau_{\dl / 4 } (d (w) )  > s 
}. 
\end{equation*} 
By combining this fact with \eqref{eq:rep2}, 
the Portmanteau theorem 
together with 
\eqref{eq:rep1}, \eqref{eq:rep2} and \eqref{eq:rep3} 
yields 
\begin{multline} \label{eq:rep4}
\limsup_{\a \to 0}
\E \cbra{ 
  \psi ( d ( \mathbf{X}^\a (s) ) )
  \; ; \; 
  \{ \hat{\tau}_\dl > s \} 
    \cap
  \{ \hat{\sg}_R > s \} 
    \cap 
  ( E_\ep^\a )^c
}
\\
\le 
\E \cbra{ 
  \psi ( \tilde{\rho} (s) + C ( \ep , t ) ) 
  \; ; \; 
  \tau_{\dl / 4} ( \tilde{\rho} ) > s
}
+ 
\P [ \sg_{R^*} ( \tilde{\rho} ) \le s ] . 
\end{multline}
In a similar way as \eqref{eq:rep2} and \eqref{eq:rep3}, 
we obtain
\begin{multline} \label{eq:rep5}
\E \cbra{ 
  \psi ( \tilde{\rho} (s) + C ( \ep , t ) ) 
  \; ; \; 
  \tau_{\dl / 4} ( \tilde{\rho} ) > s 
}
+ 
\P [ \sg_{R^*} ( \tilde{\rho} ) \le s ]  
\\
\le 
\E \cbra{ 
  \psi ( \rho (s) + C ( \ep , t ) ) 
  \; ; \; 
  \tau_{\dl / 4} ( \rho ) > s 
}
+ 
2 \P [ \sg_{R^*} ( \rho ) \le s ] .
\end{multline}
Here we used the fact $\psi \le 1$. 
Since $\mathbf{X}^\a$ converges in law to $\tilde{\mathbf{X}}$ 
in $\mathscr{C}$, 
by applying the Portmanteau theorem to \eqref{eq:loc1} 
together with \eqref{eq:rep4} and \eqref{eq:rep5}, 
we obtain 
\begin{align} \nonumber 
\E \cbra{ \psi ( d ( \mathbf{X} (s) ) ) } 
&  = 
\lim_{ \a \to 0 } 
\E \cbra{ \psi ( d ( \mathbf{X}^\a (s) ) ) }
\\ \nonumber
&  \le 
\E \cbra{ 
  \psi ( \rho (s) + C ( \ep , t ) ) 
  \; ; \; 
  \tau_{\dl / 4} ( \rho ) > s 
}
  + 
2 \P [ \sg_{R^*} ( \rho ) \le s ]
\\ \nonumber 
&  \hspace{2em} 
  + 
\limsup_{\a \to 0} \P [ \hat{\sg}_R \le s ]
  + 
\E \cbra{ 
  \psi ( d ( \mathbf{X} (s) ) )
  \; ; \; 
  \tau_\dl ( d ( \mathbf{X} (\cdot) ) ) \le s 
} 
  + 
\ep . 
\end{align}
By letting $\ep \to 0$ in this inequality, 
we obtain 
\begin{multline} \label{eq:bound1}
\E \cbra{ 
  \psi ( d ( \mathbf{X} (s) ) )
  \; ; \; 
  \tau_\dl ( d ( \mathbf{X} (\cdot) ) ) > s 
} 
\\
\le 
\E \cbra{ 
  \psi ( \rho (s) ) 
  \; ; \; 
  \tau_{\dl / 4} ( \rho ) > s 
}
  + 
2 \P [ \sg_{R^*} ( \rho ) \le s ]
  + 
\limsup_{\a \to 0} \P [ \hat{\sg}_R \le s ]. 
\end{multline}
By the definition of $\mathbf{X}^*$ and $\tau^*$, 
we have 
\begin{align} \nonumber 
\lim_{\dl \to 0} 
\E \cbra{
  \psi ( d ( \mathbf{X} (s) ) )
  \; ; \; 
  \tau_\dl ( d ( \mathbf{X} (\cdot) ) ) > s 
}   
&  = 
\lim_{\dl \to 0}
\E \cbra{
  \psi ( d ( \mathbf{X}^* (s) ) )
  \; ; \; 
  \tau_\dl ( d ( \mathbf{X}^* (\cdot) ) ) > s 
}   
\\ \nonumber 
&  = 
\E \cbra{
  \psi ( d ( \mathbf{X}^* (s) ) )
  \; ; \; 
  \tau^* > s 
}  
\\ \label{eq:Xrep}
&  = 
\E \cbra{ \psi ( d ( \mathbf{X}^* (s) ) ) } . 
\end{align}
Here the last equality follows from $\psi (0) = 0$. 
Similarly we obtain 
\begin{equation} \label{eq:Rrep}
\lim_{\dl \to 0} 
\E \cbra{ 
  \psi ( \rho (s) ) 
  \; ; \; 
  \tau_{\dl / 4} ( \rho ) > s 
}
  = 
\E \cbra{ 
  \psi ( \rho (s) ) 
  \; ; \; \tau_0 ( \rho ) > s
}. 
\end{equation}
Thus, by combining \eqref{eq:bound1} 
with \eqref{eq:Xrep} and \eqref{eq:Rrep} 
and 
by tending $R \to \infty$ with \eqref{eq:cons} in mind, 
we obtain 
\[
\E \cbra{ \psi ( d ( \mathbf{X}^* (s) ) ) }
\le 
\E \cbra{ 
  \psi ( \rho (s) ) 
  \; ; \; 
  \tau_0 ( \rho ) > s  
}. 
\]
Here we used the fact that $\rho$ cannot hit $\bar{R}$ in finite time. 
Hence the assertion holds. 
\end{Proof}

To complete the proof of Theorem~\ref{th:main}, 
we will use a regularity result on $\ph_t$ 
in the next section. 
As you will see, 
all the arguments in the next section are 
independent of this section. 
Thus there are no danger of circular arguments. 

\begin{tProof}{\protect{Theorem~\ref{th:main}}}
By virtue of Proposition~\ref{prop:concave} \ref{item:smooth} below, 
we can apply Proposition~\ref{prop:main0} 
with $\psi = \ph_t$. 
Thus we obtain 
\begin{equation} \label{eq:mono0}
\E [ \ph_t ( d ( \mathbf{X}^* (s) ) ) ] 
\le 
\E [ \ph_t ( \rho (s) ) \; ; \; \tau_0 ( \rho ) > s ]. 
\end{equation}
Since $- \rho_{K,N,a} \stackrel{d}{=} \rho_{K,N,-a}$ holds, 
a process 
$\tilde{\mathbf{\rho}}^* = ( \tilde{\rho}^{(1)} , \tilde{\rho}^{(2)} )$ 
given by 
\begin{equation*}
\tilde{\mathbf{\rho}}^* (t) 
: = 
\begin{cases} 
\displaystyle
\abra{
  \frac{\rho_{K,N,a} (t)}{2} , 
  - \frac{\rho_{K,N,a} (t)}{2}
}
& \mbox{if $t < \tau_0 ( \rho_{K,N,a} )$,}
\\
\displaystyle
\abra{
  \frac{\rho_{K,N,a} (t)}{2} , 
  \frac{\rho_{K,N,a} (t)}{2} 
} 
& 
\mbox{if $t \ge \tau_0 ( \rho_{K,N,a} )$.}
\end{cases}
\end{equation*}
is a coupling of $\rho_{K,N,a}/2$ and $\rho_{K,N,-a}/2$. 
Since the reflection map 
$x \mapsto -x$ on $( - \bar{R}/2 , \bar{R}/2 )$ 
provides a reflection structure for $\rho_{K,N,a}/2$ 
in the sense in \cite{K4}, 
$\tilde{\mathbf{\rho}}^*$ is a maximal coupling 
of $\rho_{K,N,a}/2$ and $\rho_{K,N,-a}/2$, 
and $\tau_0 ( \rho_{K,N,a}/2 ) = \tau_0 ( \rho_{K,N,a} )$ is 
the coupling time of $\tilde{\mathbf{\rho}}^*$. 
Thus Remark~\ref{rem:comparison} \ref{reflection} 
yields 
\begin{equation} \label{eq:max}
\ph_t ( | \tilde{\mathbf{\rho}}^* (s) | ) = \ph_{t+s} (a). 
\end{equation} 
Since the definition of $\tilde{\mathbf{\rho}}^*$ implies   
\[
\E [ \ph_t ( \rho (s) ) \; ; \; \tau_0 ( \rho ) > s ] 
= 
\E [ \ph_t ( | \tilde{\mathbf{\rho}}^* (s) | ) ], 
\]
the combination of it with \eqref{eq:mono0} and \eqref{eq:max} 
deduces the conclusion. 
\end{tProof}

\section{Properties of the cost function}
\label{sec:cost}

Let us define $\chi : [ 0, \infty ] \to [ 0 , 1 ]$ 
by 
\[
\chi (r) 
 := 
\frac{1}{\sqrt{2\pi}} \int_{-r}^r \e^{- u^2 / 2} du 
\]
and $\chi ( \infty ) = 1$. 
We can easily verify that 
$\chi$ is increasing and concave. 
The first goal of this section is the following expression of 
$\ph_t (a)$: 

\begin{Prop} \label{prop:express}
For each $N \in [ 2, \infty ]$, $K \in \R$ and $t \ge 0$,
there exists a probability measure $\zeta_{t,K,N}$ 
on $[0, \infty )$ such that 
\begin{equation} \label{eq:express}
\ph_t (a) : =
\int_{[0, \infty)} \chi \abra{ \frac{a}{2\sqrt{2u}} } \zeta_{t,K,N} ( du )  
\end{equation}
holds for each $a \in [ 0 , \infty )$. 
In addition, we can take $\zeta_{t,K,N}$ so that 
it is continuous in $t$ 
with respect to the topology of weak convergence. 
\end{Prop} 
The expression \eqref{eq:express} will be used to study 
some properties of $\ph_t (a)$ in Proposition~\ref{prop:concave}. 
We divide the proof of Proposition~\ref{prop:express} 
into the following two lemmata; 
Lemma~\ref{lem:Eucl} when $N = \infty$ or $K = 0$
and Lemma~\ref{lem:tc2} when $N < \infty$ and $K \neq 0$. 
We will give an expression of $\zeta_{t,K,N}$ there. 

\begin{Lem} \label{lem:Eucl}
Suppose $N=\infty$ or $K = 0$. 
Then 
\begin{align*}
\ph_t (a) 
& = 
\chi \abra{ \frac{a}{2 \sqrt{2 \eta (t)}} },  
\end{align*}
where $\eta (t) = \eta_K (t)$ is given by 
\begin{align*}
\eta_K (t) 
& := 
\begin{cases}
\displaystyle \frac{\e^{2Kt} - 1}{2K} & K \neq 0 ,
\\
\displaystyle t & K = 0 .
\end{cases}
\end{align*}
In particular, Proposition~\ref{prop:express} 
holds with $\zeta_{t, K, \infty} = \zeta_{t,0, N} = \dl_{\eta (t)}$. 
\end{Lem}

\begin{Proof}
In this case, 
$\displaystyle \rho_t : = \e^{-Kt} a +  2 \sqrt{2} \int_0^t \e^{K(s-t)} d \b_s$ 
holds. 
By the martingale representation theorem, 
$\int_0^t \e^{Ks} d \b_s$ 
can be written as a deterministic time-change of 
a standard one-dimensional Brownian motion. 
By using this fact together with the expression of 
the hitting time distribution of the Brownian motion, 
the desired expression of $\ph_t$ follows. 
\end{Proof}

To consider the case $N < \infty$, 
we begin with the following auxiliary lemma: 
\begin{Lem} \label{lem:tc1}
Suppose $N < \infty$. 
Let $\b^i (t)$ be the standard one-dimensional Brownian motion 
for $i=1,2$. 
Let $\te (t)$ be the solution to 
the following stochastic differential equation: 
\begin{align*}
d \te (t)
& = 
\sqrt{2} d \b^1 (t) 
+ \abra{ 
  \frac{N-2}{t_{K/(N-1)}( \te (t))} 
  - \frac{K}{N-1} t_{K/(N-1)} ( \te (t)) 
} dt , 
\\
\te (0) & = 0 .
\end{align*}
Let us define $\Xi (t)$ by 
\begin{equation} \label{eq:tc}
\Xi (t) : = 2 s_{K/(N-1)}^{-1} \abra{ 
  c_{K/(N-1)} \abra{ \te (t) } 
  s_{K/(N-1)} \abra{ 
    \frac{a}{2} + \sqrt{2} \int_0^t \frac{d \b^2 (s)}{c_{K/(N-1)} ( \te (s) )} 
  }
}. 
\end{equation}
Then $\Xi$ has the same law as $\rho$. 
\end{Lem}
The alternative expression of $\rho$ in the last lemma 
comes from a skew-product expression of 
the distance between two Brownian motions 
coupled by reflection on a sphere. 
For explaining a heuristic idea behind it, 
we assume $N \in \N$, $K = N-1$ and $Z \equiv 0$ for a while. 
We identify the sphere $\mathbb{S}^N$ of constant sectional curvature 1 
with an unit sphere in $\R^{N+1}$ as a submanifold. 
Let $H$ be a (uniquely determined) 2-dimensional plane in $\R^{N+1}$ 
containing origin and given starting points of the coupling 
of Brownian motions by reflection. 
Then we can decompose the Brownian motion on $\mathbb{S}^N$ 
into the ``circular part'', that is, the projection to $H$ 
and the ``complementary part'', that is, the projection to $H^\perp$. 
As a result, we can describe the distance between 
the two Brownian particles coupled by reflection 
by the scaled distance between two time-changed Brownian motions 
coupled by reflection on a circle 
whose space scaling and clock process 
are given by functionals of the complementary part. 
This description leads us to 
the expression in Lemma~\ref{lem:tc1}. 
Moreover, once we obtained this expression, 
we can verify it valid even when $N \notin \N$ or $K < 0$
as we will see in the following proof of Lemma~\ref{lem:tc1}. 
For the skew product expression of spherical Brownian motions, 
see \cite{Ito-McK,Pau-Rog_skew}, for example. 

\begin{tProof}{\protect{Lemma~\ref{lem:tc1}}}
For simplicity of notations, 
we denote $K / (N-1)$ by $\bar{K}$ in this proof. 
Let 
$\hat{\Xi} (t) : = s_{\bar{K}} ( \Xi (t) / 2 )$ 
and 
$\hat{\rho} (t) : = s_{\bar{K}} ( \rho (t) / 2)$. 
It suffices to show that 
both $\hat{\Xi} (t)$ and $\hat{\rho} (t)$ solves 
the following stochastic differential equation 
\begin{align} \label{eq:hyper0}
d z(t) 
& = 
\sqrt{ 1 - \bar{K} z(t)^2 } d w(t)  
- \frac{N\bar{K}}{2} z(t) dt 
\end{align}
for a standard Brownian motion $w(t)$. 
By the It\^{o} formula together with \eqref{eq:rho0}, 
we can easily verify that 
$\hat{\rho}$ solves \eqref{eq:hyper0} with $w(t) = \b (t)$. 
The It\^o formula together with \eqref{eq:tc} yields 
\begin{align*}
d \hat{\Xi} (t) 
& = 
- 
\bar{K} s_{\bar{K}} ( \te (t) ) 
s_{\bar{K}} \abra{ 
  \frac{a}{2} 
  + 
  \sqrt{2} \int_0^t \frac{d \b^2 (s)}{c_{\bar{K}} ( \te (s) )}  
} 
d \te (t) 
\\
& \hspace{2em}
+ \sqrt{2}
c_{\bar{K}} \abra{ 
  \frac{a}{2} 
  + 
  \sqrt{2} \int_0^t \frac{d\b^2 (s)}{c_{\bar{K}} ( \te (s) )} 
} d\b^2 (t)
\\
& \qquad 
- 
\bar{K} \abra{ 
  c_{\bar{K}} ( \te (t) ) 
  + 
  \frac{1}{c_{\bar{K}} ( \te (t) )}
}
s_{\bar{K}} \abra{ 
  \frac{a}{2} 
  + \sqrt{2} \int_0^t \frac{d \b^2 (s)}{c_{\bar{K}} ( \te (s) )}  
} 
dt 
\\
& = 
- \sqrt{2}
\bar{K} s_{\bar{K}} ( \te (t) ) 
s_{\bar{K}} \abra{ 
  \frac{a}{2} 
  + 
  \sqrt{2} \int_0^t \frac{d \b^2 (s)}{c_{\bar{K}} ( \te (s) )} 
} 
d \b^1 (t) 
\\
& \qquad 
+ \sqrt{2} 
c_{\bar{K}} \abra{ 
  \frac{a}{2} 
  + 
  \sqrt{2} \int_0^t \frac{d\b^2 (s)}{c_{\bar{K}} ( \te (s) )} 
} d\b^2 (t)
- N\bar{K} 
\hat{\Xi} (t) dt . 
\end{align*} 
Here we used the relation 
$c_{\bar{K}} (r)^2 + \bar{K} s_{\bar{K}} (r)^2 = 1$, 
which holds for any $K \in \R$, 
to obtain the last equality. 
By a direct computation, we have 
\begin{align*}
& \abra{ 
  \bar{K} s_{\bar{K}} ( \te (t) ) 
  s_{\bar{K}} \abra{ 
    \frac{a}{2} 
    + 
    \sqrt{2} \int_0^t \frac{d \b^2 (s)}{ c_{\bar{K}} ( \te (s) ) } 
  }
}^2 
+ 
c_{\bar{K}} \abra{ 
  \frac{a}{2} 
  + 
  \sqrt{2} \int_0^t \frac{d \b^2 (s)}{c_{\bar{K}} ( \te (s) )} 
}^2 
= 1 - \bar{K} \hat{\Xi} (t)^2 . 
\end{align*}
Note that $1 - \bar{K} \hat{\Xi} (t)^2 > 0$ holds 
for any $t \ge 0$ almost surely 
since $\te (t)$ never hits $\bar{R} / 2$. 
Thus, $\hat{\Xi}$ solves \eqref{eq:hyper0} with $w (t) = \b^* (t)$ 
given by 
\begin{multline*}
\b^* (t) := \int_0^t 
\frac{1}{\sqrt{1 - \bar{K} \hat{\Xi} (t)^2}}
\Bigg( 
  - 
  \bar{K} s_{\bar{K}} ( \te (s) ) 
  s_{\bar{K}} \abra{ 
    \frac{a}{2} 
    + 
    \sqrt{2} \int_0^s \frac{d \b^2 (u)}{c_{\bar{K}} ( \te (u) )}  
  } 
  d \b^1 (s) 
\\
  + c_{\bar{K}} \abra{ 
    \frac{a}{2} 
    + 
    \sqrt{2} \int_0^s \frac{d\b^2 (u)}{c_{\bar{K}} ( \te (u) )} 
  } d\b^2 (s)
\Bigg)   
\end{multline*}
and hence the conclusion follows. 
\end{tProof}

\begin{Lem} \label{lem:tc2}
Suppose $N < \infty$. 
We denote the law of 
$\int_0^t c_{K/(N-1)} ( \te (s))^{-2} ds$ 
by $\zeta_{t,K,N}$ for each $t \ge 0$, 
where $\te (\cdot)$ is as in Lemma~\ref{lem:tc1}. 
Then the conclusion of 
Proposition~\ref{prop:express} holds true. 
\end{Lem}

\begin{Proof}
The continuity in $t$ of $\zeta_{t,K,N}$ directly follows from the definition. 
Let $\Xi$ be as in Lemma~\ref{lem:tc1}. 
By the martingale representation theorem, 
there exists a standard one-dimensional Brownian motion $B(t)$ 
such that 
\[
B \abra{ 2 \int_0^\cdot \frac{ds}{ c_{K/(N-1)} (\te (s))^2 } } 
\stackrel{d}{=} 
\sqrt{2} \int_0^\cdot \frac{d\b^2 (s)}{c_{K/(N-1)} (\te (s))} 
\]
holds. 
Since $\b^1$ and $\b^2$ are independent, 
$B (\cdot)$ behaves as a standard Brownian motion
even under the conditional probability measure 
$\P[ \;\cdot \; | \sg ( \b^1 ) ]$. 
Thus the definition of $\ph_t (a)$ and Lemma~\ref{lem:tc1} yield 
\begin{align*}
\ph_t (a) 
& = 
\P \cbra{ \inf_{0 \le s \le t} \rho (s) > 0 } 
= 
\P \cbra{ 
  \inf_{0 \le s \le t} 
  \abra{ 
    \frac{a}{2}  
    + 
    \sqrt{2} \int_0^s \frac{d \b^2 (u)}{ c_{K/(N-1)} (\te (u))} 
  }
  > 0 
}
\\
& = 
\P \cbra{ 
  \inf 
  \bbra{
    \frac{a}{2} + B(s) \; \left| \; 
    0 \le s \le 2 \int_0^t \frac{du}{ c_{K/(N-1)} (\te (u))^2} 
    \right.
  } 
  > 0
}
\\
& = 
\E \cbra{ 
  \P \cbra{ \left.
    \inf 
    \bbra{
      \frac{a}{2} + B(s) \; \left| \; 
          0 \le s \le 2 \int_0^t \frac{du}{ c_{K/(N-1)} (\te (u))^2} 
      \right.
    } 
    > 0
    \right| \sg ( \b^1 ) 
  }
}
\\
& = 
\E \cbra{ 
  \chi \abra{  
    \frac{a}{2\sqrt{2}} 
    \abra{ 
      \int_0^t \frac{du}{ c_{K/(N-1)} (\te (u))^2} 
    }^{-1/2}
  }
}. 
\end{align*}
Hence the desired result holds. 
\end{Proof}

Now we state some  
consequences of the expressions of $\ph_t (a)$ 
in Proposition~\ref{prop:express}: 

\begin{Prop} \label{prop:concave}
\begin{enumerate}
\item \label{item:t-conti}
For $a \in [ 0 , \bar{R} ]$, 
$[ 0, \infty ) \ni t \mapsto \ph_t (a)$ is continuous. 
\item \label{item:smooth}
$\ph_t$ is continuous on $\overline{[ 0 , \bar{R} )}$ 
and smooth on $( 0 , \bar{R} )$ for $t > 0$.
\item \label{item:concave}
$\ph_t$ is concave on $\overline{[ 0, \bar{R} )}$ for $t \ge 0$. 
\item \label{item:dim-compare}
For $t \ge 0$, 
$K, K' \in \R$ with $K \ge K'$, 
$N, N' \in [ 2, \infty ]$ with $N \le N'$ 
and $a \in [ 0 , \bar{R}_{K,N} ]$, 
\begin{equation*}
\ph^{K,N}_t (a) \le \ph^{K',N'}_t (a).
\end{equation*} 
\item \label{item:diff_at_zero}
For $t > 0$, $K \in \R$ and $N \in [ 2 , \infty ]$, 
$\ph^{K,N}_t$ is differentiable at $0$. 
Moreover, for $K' \in \R$ and $N' \in [2 , \infty ]$
with $K' \le K$, $N' \ge N$, 
\begin{equation} \label{eq:diff-compare}
( \ph^{K,N}_t )' (0) 
= 
\int_{[0,\infty)} \frac{\zeta_{t,K,N} (du)}{4 \sqrt{\pi u}} 
\le 
( \ph^{K',N'}_t )' (0). 
\end{equation}
In particular, 
$( \ph^{K,N}_t )' (0) \le ( \ph^{K, \infty} )' (0) = ( \pi \eta(t) )^{-1/2}/4$. 
Here $\eta (t) = \eta_K (t)$ is as in Lemma~\ref{lem:Eucl}. 
\item \label{item:t-short}
For $K \in \R$ and $N \in [ 2 , \infty ]$, 
$\displaystyle 
\lim_{t \downarrow 0} 
\sqrt{t} ( \ph_t^{K,N})' (0) 
= 
\frac{1}{4\sqrt{\pi}}
$. 
\end{enumerate}
\end{Prop}

\begin{Proof}
\ref{item:t-conti}
It is obvious by the continuity of $\zeta_{t,K,N}$ in \eqref{eq:express}. 

\ref{item:smooth}
Note that the derivative of $\chi ( a / ( 2 \sqrt{2u} ) )$ of any order 
with respect to $a$-variable is a bounded function of $u$ 
for $a \in ( 0, \bar{R} )$ in \eqref{eq:express}. 
Thus the dominated convergence theorem 
yields that $\ph_t$ is smooth on $( 0 , \bar{R} )$. 
We can show 
the continuity of $\ph_t$ on $\overline{[0, \bar{R} )}$ 
similarly. 

\ref{item:concave}
Since $\ph_0 = 1_{ ( 0 , \infty ) }$ by definition, 
it is obviously concave. 
Thus it suffices to consider the case $t > 0$. 
As we did in the proof of \ref{item:smooth}, 
we can compute $\ph_t'' (a)$ at $a \in ( 0 , \bar{R} )$ 
by using the dominated convergence theorem.  
Since $\chi$ is concave, $\ph_t'' (a) \le 0$ 
and hence the conclusion holds because $\ph_t$ is continuous 
on $\overline{[ 0 , \bar{R} )}$ by \ref{item:smooth}. 

\ref{item:dim-compare}
By a direct computation, 
we can verify $\Psi_{K,N} (u) \le \Psi_{K',N'} (u)$ 
for any $u \in [ 0 , \bar{R}_{K,N} )$. 
Thus the comparison theorem for stochastic differential equations 
(see \cite{Ik-Wat} for instance) 
yields that $\rho_{K,N,a} (t) \le \rho_{K',N',a'} (t)$ 
for $a' > a$ and $t \ge 0$. 
It implies $\ph^{K,N}_t (a) \le \ph^{K',N'}_t (a')$ 
by the definition of $\ph^{K,N}_t (a)$. 
Since $\ph_t$ is continuous, the asserted inequality follows 
by tending $a' \downarrow a$. 

\ref{item:diff_at_zero}
Since $\chi$ is concave and $\chi (0) = 0$, 
$\chi(r) / r$ is nonincreasing. 
Thus the monotone convergence theorem yields 
\begin{equation*}
\lim_{a \downarrow 0} \frac{\ph^{K,N}_t (a) - \ph^{K,N}_t (0)}{a}
= 
\lim_{a \downarrow 0} \frac{\ph^{K,N}_t (a)}{a} 
= 
\int_{[0,\infty)} \frac{\zeta_{t,K,N} (du)}{4 \sqrt{\pi u}}. 
\end{equation*}
By combining this identity with \ref{item:dim-compare}, we obtain 
\begin{equation*}
\int_{[0,\infty)} \frac{\zeta_{t,K,N} (du)}{4 \sqrt{\pi u}} 
\le 
\int_{[0,\infty)} \frac{\zeta_{t,K',N'} (du)}{4 \sqrt{\pi u}}
\le 
\int_{[0,\infty)} \frac{\zeta_{t,K',\infty} (du)}{4 \sqrt{\pi u}}
= 
\frac{1}{4 \sqrt{\pi \eta_{K'} (t)}} < \infty 
\end{equation*}
and hence the conclusion follows.  

\ref{item:t-short} 
We use the expression of $( \ph^{K,N}_t )' (0)$ in \ref{item:diff_at_zero}. 
When $N = \infty$, it easily follows from Lemma~\ref{lem:Eucl}. 
Next we consider the case $K > 0$ with the expression of $\zeta_{t,K,N}$ 
given in Lemma~\ref{lem:tc2}. 
By the definition of $\te (t)$ in Lemma~\ref{lem:tc1}, 
$c_{K/(N-1)} ( \te (t) ) \in ( 0 , 1 ]$ holds for each $t > 0$. 
Thus we have $\int_0^t ( c_{K/(N-1)} ( \te (s) ) )^{-2} ds \ge t$ 
and therefore the dominated convergence theorem yields 
\begin{align*}
\lim_{t \downarrow 0} \sqrt{t} ( \ph^{K,N}_t )' (0) 
& = 
\lim_{t \downarrow 0}
\frac{1}{4 \sqrt{\pi}} 
\E \cbra{ 
  \abra{ 
    \frac{1}{t} \int_0^t \frac{ds}{c_{K/(N-1)} ( \te (s))^2} 
  }^{-1/2}
}
= 
\frac{1}{4 \sqrt{\pi}}. 
\end{align*}
Finally, for the general $K \in \R$ and $N \in [ 2 , \infty )$, 
the conclusion follows from \eqref{eq:diff-compare} 
together with the above-mentioned two cases. 
\end{Proof}

Since $\ph_r (0) = 0$, 
Proposition~\ref{prop:concave} \ref{item:concave} 
yields the following corollary:  

\begin{Cor} \label{cor:distance}
We have $\ph_r (a + a') \le \ph_r (a) + \ph_r (a')$ 
for $r > 0$ and $a, a' \ge 0$. 
In particular, for $t > 0$, 
$\ph_t ( d ( \cdot , \cdot ) )$ is a bounded distance function 
being compatible with the topology on $M$. 
\end{Cor}

Though the preparation of the proof of Theorem~\ref{th:main2}
is already finished in Proposition~\ref{prop:concave}, 
we will discuss further properties of $\ph_t$ 
in the rest of this section. 
First we will study more explicit expression of $\ph_t (a)$ 
than the one in Lemma~\ref{lem:tc2} 
in the case $N < \infty$ and $K \neq 0$. 
Lemma~\ref{lem:hyper1} and Corollary~\ref{cor:hyper2} 
below study the case $K < 0$. 
Based on the expression of the Brownian motion on 
the hyperbolic space 
by a stochastic differential equation 
(see \cite{Matsu-Yor_II}, for instance), 
we can show the following 
in a similar way as Lemma~\ref{lem:tc1}: 

\begin{Lem} \label{lem:hyper1}
Suppose $N < \infty$ and $K < 0$. 
Let $\b^1 (t)$ and $\b^2 (t)$ be independent, 
one-dimensional standard Brownian motions. 
Let $\Xi' (t)$ and $\te' (t)$ be given by 
\begin{align*}
\te' (t) 
& : = 
\exp \abra{ 
  \sqrt{\frac{-2 K}{N-1}} \b^1 (t) 
  + Kt 
}, 
\\
\Xi' (t) 
& : = 
s_{K/(N-1)} \abra{ \frac{a}{2} } 
+ 
\sqrt{2} \int_0^t \te' (s) d \b^2 (s) . 
\end{align*}
Then $2 s_{K/(N-1)}^{-1} ( \Xi' (\cdot) / \te' (\cdot) )$ 
has the same law as $\rho$. 
\end{Lem}
\begin{Proof}
As in the proof of Lemma~\ref{lem:tc1}, 
we denote $K/ (N-1)$ by $\bar{K}$. 
We already know in the proof of Lemma~\ref{lem:tc1} 
that $s_{\bar{K}} (\rho (t) / 2 )$ 
solves the stochastic differential equation \eqref{eq:hyper0} 
with $w (t) = \b (t)$. 
Thus it suffices to show that 
$\Xi' (t) / \te' (t)$ also solves \eqref{eq:hyper0} 
for a standard Brownian motion $w (t)$. 
The It\^o formula yields 
\begin{align*}
d \abra{ \frac{\Xi' (t)}{\te' (t)} }
& = 
- \frac{\Xi' (t)}{\te' (t)^2} d\te' (t) 
+ \frac{\Xi' (t)}{\te' (t)^3} d \dbra{ \te' }(t) 
+ \frac{1}{\te' (t)} d \Xi' (t)
\\
& = 
- \sqrt{-\bar{K}} \frac{\Xi' (t)}{\te' (t)} 
  \abra{
    \sqrt{2} d\b^1 (t) 
    - (N-2)\sqrt{-\bar{K}} dt 
  } 
- 2 \bar{K} \frac{\Xi' (t)}{\te' (t)} dt 
+ \sqrt{2} d\b^2 (t) 
\\
& = 
\sqrt{2} d\b^2 (t) 
- \sqrt{-2 \bar{K}}\frac{\Xi' (t)}{\te' (t)} d \b^1 (t) 
- N\bar{K} \frac{\Xi' (t)}{\te' (t)} dt . 
\end{align*}
Thus $\Xi' (t) / \te' (t)$ solves \eqref{eq:hyper0} 
with $w(t) = \b^{**} (t)$ given by 
\[
\b^{**} (t) 
:= 
\int_0^t 
\abra{ 
  1  - \bar{K} \abra{ \frac{ \Xi' (s) }{ \te' (s) } }^2 
}^{-1/2} 
\abra{
  d\b^2 (s) 
  - 
  \sqrt{-\bar{K}} \abra{ \frac{ \Xi' (s) }{ \te' (s) } } 
  d\b^1 (s) 
}
. 
\]
\end{Proof}

\begin{Cor} \label{cor:hyper2}
Suppose $N < \infty$ and $K < 0$. Then 
\begin{align} \label{eq:pre-hyper}
\ph_t (a) = 
\E \cbra{ 
  \chi \abra{ 
    \frac{1}{2\sqrt{2}} s_{K/(N-1)} \abra{ \frac{a}{2} }
    \abra{
      \int_0^t \te' (s)^2 ds  
    }^{-1/2}
  }
}, 
\end{align}
where $\te' (t)$ is as in Lemma~\ref{lem:hyper1}. 
Moreover, 
\begin{align*}
\ph_t (a) 
& = 
\int_{- \infty}^\infty \int_0^\infty 
\chi \abra{
  \frac12 \sqrt{ \frac{-K}{(N-1)u} } 
  s_{K/(N-1)} \abra{ \frac{a}{2} }
} 
\\
& \hspace{4em}
\times 
\exp \abra{ 
  \frac{(N-1)}{2} ( Kt - x ) 
  - 
  \frac{ 1 + \e^{2x}}{2u}
}
\vartheta \abra{ \frac{\e^x}{u} , \frac{-2 K t}{N-1} } 
\frac{d u}{u} dx , 
\end{align*}
where 
\begin{equation*}
\vartheta ( r , t ) : = 
\frac{r}{2 \pi^3 t}
\e^{\pi^2/(2t)} 
\int_0^\infty 
\e^{ - \xi^2 / (2t)} 
\e^{-r \cosh (\xi)}
\sinh (\xi) 
\sin \abra{ \frac{ \pi \xi }{t} }
d \xi .  
\end{equation*}
\end{Cor}
\begin{Proof}
Let $\Xi' (t)$ and $\te' (t)$ be as in Lemma~\ref{lem:hyper1}. 
By the martingale representation theorem, 
there exists a Brownian motion $B(t)$ 
such that 
\[
\Xi' (t) 
\stackrel{d}{=} 
s_{K/(N-1)} \abra{ \frac{a}{2} } 
+ B \abra{ 2 \int_0^t \te' (s)^2 ds }.
\]
Hence, as in the proof of Lemma~\ref{lem:tc2}, 
the definition of $\ph_t (a)$ and Lemma~\ref{lem:hyper1} yield 
\begin{align*}
\ph_t (a) 
& = 
\P \cbra{ \inf_{0 \le s \le t} \rho (s) > 0 } 
= 
\P \cbra{ \inf_{0 \le s \le t} \Xi' (s) > 0 }
\\
& = 
\P \cbra{ 
  \inf 
  \bbra{
    s_{K/(N-1)} \abra{ \frac{a}{2} }  + B(s) \; \left| \; 
    0 \le s \le 2 \int_0^t \te' (u)^2 du 
    \right.
  } 
  > 0
}
\\
& = 
\E \cbra{ 
  \chi \abra{  
    \frac{1}{2\sqrt{2}} s_{K/(N-1)} \abra{ \frac{a}{2} }
    \abra{ 
      \int_0^t \te' (u)^2 du 
    }^{-1/2}
  }
}. 
\end{align*}
This is nothing but \eqref{eq:pre-hyper}. 
Now the conclusion follows 
by using an explicit expression of 
the distribution of $\int_0^t \te' (u)^2 du$ 
in \cite[Theorem~4.1]{Matsu-Yor_I} 
(also see references therein). 
\end{Proof}

In the case $K > 0$, we use several properties on 
the Gegenbauer, or ultraspherical, polynomials 
to obtain alternative expression of $\ph_t$ 
in Lemma~\ref{lem:Gegenbauer} below. 
We refer to \cite{Szego_OP} 
for basics on Gegenbauer polynomials. 

\begin{Lem} \label{lem:Gegenbauer}
Suppose $N < \infty$ and $K > 0$. Then, 
for all $a \in [ 0, \bar{R} ]$, 
\[
\ph_t ( a ) 
= 
\sum_{n=0}^\infty 
\e^{ - (2n+1)(2n+N) K t / (N-1) } 
\frac{ (-1)^n (4n + N + 1) }{ \pi (2n+N) } 
B \abra{ \frac{N-1}{2} , n + \frac12 }
P_{2n+1} (\tilde{a} ), 
\]
where $B ( \cdot, \cdot )$ is the Beta function, 
$
\tilde{a} 
 : = 
\sin ( \sqrt{K/(N-1)} a /2 ) 
$
and 
$P_n (x)$ is the $n$-th Gegenbauer polynomial of parameter $(N-1)/2$. 
\end{Lem}

\begin{Proof}
Let us define $\hat{\rho} (t)$ by 
\[
\hat{\rho} (t) 
:= 
\sin \abra{ 
  \frac12 \sqrt{\frac{K}{N-1}}
  \rho \abra{ \frac{(N-1)t}{2 K} }
} . 
\]
Then $\hat{\rho} (t)$ solves the following stochastic differential equation: 
\begin{align*}
d \hat{\rho} (t) 
& = 
\sqrt{ 1- \hat{\rho} (t)^2 } d \b (t) 
- \frac{N}{2} \hat{\rho}(t) dt , 
\\
\hat{\rho} (0) 
& = 
\tilde{a}, 
\end{align*} 
where $\b (t)$ is a one-dimensional standard Brownian motion. 
Thus $\hat{\rho}(t)$ is the Legendre process, 
or the diffusion process on $(-1, 1)$ generated 
by $L_N$ given as follows:  
\[
L_N := \frac12 (1-x^2) \dfrac{\partial^2 }{\partial x^2} 
- \frac{N}{2} x \dfrac{\partial }{\partial x} .
\]  
It is well-known that $\mu_N (dx) = ( 1 - x^2 )^{N/2 - 1} dx$ is 
the symmetrizing measure of $L_N$. Moreover, 
$L_N$ is essentially selfadjoint on $L^2 (\mu_N)$, 
the spectra of $L_N$ on $L^2 ( \mu_N )$ is 
$\{ - n ( n + N - 1 ) / 2 \}_{n \in \N_0}$ 
all of which are eigenvalues of multiplicity one, and 
the normalized eigenfunction corresponding 
to the $n$-th eigenvalue is $n$-th normalized Gegenbauer polynomials
$\bar{P}_n$ defined by $\bar{P}_n (x) = Z_n^{-1} P_n (x)$ and 
\begin{align*}
Z_n 
& : = 
\bbra{ \int_0^1 P_n (x)^2 \mu_N (dx) }^{1/2} 
= 
\frac{ 
  2^{1- N/2} \sqrt{\pi \Gm ( n + N - 1 )}
}
{
  \sqrt{n! ( n + (N-1)/2 )} \Gm ((N-1)/2)
} .  
\end{align*} 
As a result, the transition density $p_1 ( t, x, y )$
of $\hat{\rho}(t)$ with respect to $\mu_N$ 
is given by 
\begin{equation} \label{eq:SD}
p_1 ( t , x , y ) 
= 
\sum_{n=0}^\infty 
\e^{ -n ( n + N - 1 ) t / 2 } \bar{P}_n (x) \bar{P}_n (y), 
\end{equation}
where the sum converges in $L^2 ( \mu_N \otimes \mu_N )$.  
We claim that the infinite sum 
in the right hand side of \eqref{eq:SD} 
converges uniformly in $x$ and $y$. 
Recall that the Gegenbauer polynomial $P_n$ satisfies the 
following recursion relation: 
\begin{align} \label{eq:recur}
P_n (x) 
& = 
\frac{2n + N - 3}{n} x P_{n-1} (x) 
- 
\frac{( n + N - 3 )}{n} P_{n-2} (x), 
\\ \nonumber
P_0 (x) & = 1 , \quad P_1 (x) = ( N - 1 ) x.  
\end{align}
By induction, we can easily show that 
there exist $c_0 > 0$ and $q > 0$ such that 
\begin{equation} \label{eq:G-exp}
\sup_{x \in ( -1 , 1 )} | P_n (x) | \le c_0 q^n . 
\end{equation}
(for instance, we can dominate the left hand side by 
$(N-1) 4^n \prod_{k=1}^n ( 1 + | N - 3 | / k )$). 
It is not difficult to see that there exists $c_1 > 0$ 
such that 
$Z_n^{-1} \le c_1 \sqrt{n}$ for all $n \in \N$ 
since $N \ge m \ge 2$. 
Then these estimates imply the claim. 

Now the reflection principle yields 
\begin{align*}
\ph_t (a) 
& = 
\P \cbra{ \inf_{0 \le s \le t} \rho (s) > 0 } 
\\
& = 
\P \cbra{ \inf_{0 \le s \le 2Kt/(N-1)} \hat{\rho} (s) > 0 } 
\\
& = 
\int_0^1 
\abra{ 
  p_1 ( \frac{2Kt}{N-1} , \tilde{a} , x ) 
  - 
  p_1 ( \frac{2Kt}{N-1} , - \tilde{a} , x ) 
} 
\mu_N ( dx )
\\
& = 
\sum_{n=0}^\infty 
\int_0^1 
\e^{ -n ( n + N - 1 ) Kt / (N-1) } 
\abra{ 
  \bar{P}_n (\tilde{a}) - \bar{P}_n ( - \tilde{a}) 
} \bar{P}_n (x) \mu_N (dx)
\\
& = 
2 \sum_{n=0}^\infty 
\e^{ - ( 2n + 1 )( 2n + N ) K t / (N-1) }
P_{2n+1}^* ( \tilde{a} )
\int_0^1 P_{2n+1}^* (x) \mu_N ( dx ) .
\end{align*}
The Rodrigues formula for the Gegenbauer polynomial asserts 
\begin{equation*}
P_{n} (x) 
= 
\frac{( - 2 )^{n}}{n!} 
\frac{ 
  \Gm ( n + (N-1)/2 ) \Gm ( n + N - 1)
}{
  \Gm ( (N-1)/2 ) \Gm ( 2 n + N - 1 )
}
( 1 - x^2 )^{1-N/2}
\frac{d^{n}}{d x^{n}} ( 1 - x^2 )^{n + N/2 - 1}. 
\end{equation*}
By using this formula twice, we obtain 
\begin{align*}
\int_0^1 P_{2n+1} (x) \mu_N (dx) 
& = 
\frac{( - 2 )^{2n+1}}{(2n+1)!} 
\frac{ 
  \Gm ( 2n + (N+1)/2 ) \Gm ( 2n + N )
}{
  \Gm ( (N-1)/2 ) \Gm ( 4 n + N + 1 )
}
\cbra{ 
  \frac{d^{2n}}{d x^{2n}} ( 1 - x^2 )^{2n + N/2} 
}_{x=0}^1 
\\
& = 
\frac{(N-1)}{(2n+1)(2n + N)}
P_{2n}^{(N+2)} (0), 
\end{align*}
where $P_{2n}^{(N+2)}$ is the $(2n)$-th Gegenbauer polynomial 
of parameter $(N+1)/2$ (associated with $L_{N+2}$). 
By the recursion formula \eqref{eq:recur}, 
we obtain 
\begin{equation*}
P^{(N+2)}_{2n} (0) 
= 
(-1)^n \frac{ \Gm ( n + ( N + 1 )/2 )}{\Gm ( (N+1)/2 ) n!}
.
\end{equation*} 
Thus the duplication formula of the Gamma function yields
\begin{align} \nonumber
\ph_t (a) 
& = 
\sum_{n=0}^\infty 
\e^{ - (2n+1)(2n+N) K t / (N-1) } 
\frac{(-1)^n 2 (N-1) \Gm ( n + ( N + 1 )/2 )}
{Z_{2n+1}^2 (2n+1)(2n + N) \Gm ( ( N + 1 )/2 ) n!}
P_{2n+1} (\tilde{a} )
\\ \nonumber 
& = 
\sum_{n=0}^\infty 
\e^{ - (2n+1)(2n+N) K t / (N-1) } 
\frac{ (-1)^n (4n + N + 1) }{ \pi (2n+N) } 
B \abra{ \frac{N-1}{2} , n + \frac12 }
P_{2n+1} (\tilde{a} ) .
\end{align}
This is nothing but the desired identity. 
\end{Proof}

Based on expressions of $\ph_t (a)$ 
in Lemma~\ref{lem:Eucl}, 
Corollary~\ref{cor:hyper2} and 
Lemma~\ref{lem:Gegenbauer}, 
we will obtain the asymptotic behavior of 
$\ph_t (a)$ as $t \to \infty$ in the following corollary: 

\begin{Cor} \label{cor:t-large}
The following convergence holds compact uniformly in $a \in [ 0, \bar{R} )$:   
\begin{enumerate}
\item \label{item:infin_nonneg}
When $N = \infty$ and $K \ge 0$, or $N < \infty$ and $K = 0$, 
\[
\lim_{t \to \infty} \sqrt{\eta_K (t)} \ph_t (a) = \frac{a}{4\sqrt{\pi}}. 
\]
In addition, this is an increasing limit. 
\item \label{item:infin_neg}
When $N = \infty$ and $K < 0$, 
\[
\lim_{t \to \infty} \ph_t (a) = \chi \abra{ \frac{a \sqrt{-K}}{2} }. 
\]
\item \label{item:fin_pos}
When $N < \infty$ and $K > 0$, 
\[
\lim_{t \to \infty} \e^{NKt/(N-1)} \ph_t (a) 
= 
\frac{(N^2 -1)}{ \pi N} 
B \abra{ \frac{N-1}{2} , \frac12 } 
\sin \abra{ \frac12 \sqrt{\frac{K}{N-1}} a}. 
\]
In addition, 
$
\sup \{ 
  \e^{NKt/(N-1)} \ph_t (a) 
  \; | \; 
  t \ge 1 , a \in [ 0 , \bar{R} ] 
\} < \infty
$. 
\item \label{item:fin_neg}
When $N < \infty$ and $K < 0$, 
\begin{equation} \label{eq:lim-hyper0}
\lim_{t \to \infty} \ph_t (a) 
 = 
\frac{1}{\Gamma ( (N-1)/2 )}
\int_0^\infty 
\chi \abra{ 
  \sqrt{\frac{-Ku}{2(N-1)}} 
  s_{K/(N-1)} \abra{\frac{a}{2}}
}
u^{(N-3)/2} \e^{-u} du .
\end{equation}
\end{enumerate}
\end{Cor}

\begin{Proof}
\ref{item:infin_nonneg} \ref{item:infin_neg} 
The convergence easily follows by elementary calculus. 
The monotonicity in $t$ in \ref{item:infin_nonneg} follows 
from the concavity of $\chi$. 
In both cases, 
the Dini theorem ensures the uniformity of the convergence 
on each compact set. 

\ref{item:fin_pos} 
The computation of the limit as well as 
the uniformity on $[ 0 , \bar{R} ]$ and 
the finiteness of the supremum 
directly follows 
from Proposition~\ref{lem:Gegenbauer} and \eqref{eq:G-exp}. 

\ref{item:fin_neg}
By \eqref{eq:pre-hyper} and the monotone convergence theorem, 
\begin{equation} \label{eq:lim-hyper1}
\lim_{t \to \infty} \ph_t (a) 
= 
\E \cbra{ 
  \chi \abra{ 
    \frac{1}{2\sqrt{2}} s_{K/(N-1)} \abra{ \frac{a}{2} }
    \abra{
      \int_0^\infty \te' (s)^2 ds  
    }^{-1/2}
  }
}. 
\end{equation}
Then the distribution $\int_0^\infty \te' (s)^2 ds$ can be described 
with the aid of \cite[Theorem~6.2]{Matsu-Yor_I} 
(also see references therein) 
to obtain \eqref{eq:lim-hyper0}. 
Since the convergence in \eqref{eq:lim-hyper1} is monotone, 
the compact uniformity of the convergence follows from 
the Dini theorem. 
\end{Proof}
 
\section{Monotonicity of transportation costs}
\label{sec:mono-cost}

Based on 
Proposition~\ref{prop:concave} 
and 
Corollary~\ref{cor:distance}, 
we will show some continuity properties 
for $\ph_t (a)$ and $\mathcal{T}_{\ph_t (d)}$ with respect to $t$ 
in the following two lemmata. 

\begin{Lem} \label{lem:time-Wconti}
Let $c_n : M \times M \to [ 0 , \infty )$ be a family of 
continuous functions converging to $c : M \times M \to [ 0 , \infty )$
pointwisely. 
Let $\mu , \nu \in \mathcal{P} (M)$.  
\begin{enumerate}
\item \label{item:W-lsup}
If $\sup_{n,x,y} c_n (x,y) < \infty$ or $c_n$ is nondecreasing in $n$, 
then 
\[
\limsup_{n \to \infty} \mathcal{T}_{c_n} ( \mu , \nu ) \le \mathcal{T}_c ( \mu , \nu ). 
\]
\item \label{item:W-linf}
If the convergence $c_n \to c$ is uniform on each compact set 
or $c_n$ is nondecreasing in $n$,  
then 
\[
\liminf_{n \to \infty} \mathcal{T}_{c_n} ( \mu , \nu ) \ge \mathcal{T}_c ( \mu , \nu ).
\] 
\end{enumerate}
\end{Lem}

\begin{Proof}
\ref{item:W-lsup}
Under the assumption on $c_n$, 
for each $\pi \in \Pi ( \mu , \nu )$, 
\[
\limsup_{n \to \infty} \mathcal{T}_{c_n} ( \mu , \nu ) 
\le 
\limsup_{n \to \infty} \int_{M \times M} c_n \, d \pi 
= 
\int_{M \times M} c \, d \pi . 
\]
Thus the assertion holds 
by taking infimum over $\pi \in \Pi ( \mu , \nu )$. 

\ref{item:W-linf}
Take a subsequence $( c_{n_k} )_k$ of $( c_n )_n$ 
so that 
\[
\lim_{k \to \infty} \mathcal{T}_{c_{n_k}} ( \mu , \nu ) 
= 
\liminf_{n \to \infty} \mathcal{T}_{c_n} ( \mu , \nu ) .
\]
Since $\Pi ( \mu , \nu )$ is compact and 
$c_n$ is continuous and nonnegative, 
a usual variational argument implies 
that there is a minimizer 
of $\mathcal{T}_{c_{n_k}} ( \mu , \nu )$, i.e.~there exists 
$\pi_k \in \Pi ( \mu , \nu )$ such that 
$\mathcal{T}_{c_{n_k}} ( \mu , \nu ) = \int_{M \times M} c_{n_k} \, d\pi_k$. 
We may assume that $\pi_k$ converges as $k \to \infty$ 
by taking a subsequence if necessary. 
We denote the limit by $\pi_\infty$. 

First we consider the case that 
$c_n$ converges to $c$ compact uniformly. 
Take $\ep > 0$ and choose a compact set $K \subset M \times M$ 
such that $\pi_k ( K^c ) < \ep$. 
Then, for any $R > 0$, the assumption on $c_n$ implies 
\begin{align*}
\lim_{k \to \infty} \mathcal{T}_{c_{n_k}} ( \mu , \nu ) 
& \ge 
\liminf_{k \to \infty} 
\int_K c_{n_k} \wg R \, d \pi_k 
\\
& \ge 
\liminf_{k \to \infty} 
\int_K c \wg R \, d \pi_k - \ep
\\
& \ge 
\liminf_{k \to \infty} \int_M c \wg R \, d \pi_k - (R+1) \ep. 
\\
& = 
\int_M c \wg R \, d \pi_\infty - (R+1) \ep. 
\end{align*} 
Thus, by taking $\ep \downarrow 0$ and $R \uparrow \infty$, 
we obtain 
\[
\lim_{k \to \infty} \mathcal{T}_{c_{n_k}} ( \mu , \nu ) 
\ge 
\int_M c \, d \pi_\infty 
\ge 
\mathcal{T}_c ( \mu , \nu ) 
\]
and hence the assertion holds. 

Next we consider the case that 
$c_n$ is nondecreasing in $n$. 
Then, for $k \in \N$, 
\[
\int_{M \times M} c_{n_k} \, d \pi_\infty 
\le  
\liminf_{l \to \infty} \int_{M \times M} c_{n_k} \, d \pi_l
\le 
\liminf_{l \to \infty} \int_{M \times M} c_{n_l} \, d \pi_l 
= 
\liminf_{n \to \infty} \mathcal{T}_{c_{n}} ( \mu , \nu ). 
\]
By taking $k \to \infty$, the monotone convergence theorem 
implies that 
\[
\mathcal{T}_c ( \mu, \nu ) 
\le 
\limsup_{k \to \infty} 
\int_{M \times M} c_{n_k} \, d \pi_\infty .  
\]
Thus, the conclusion follows by combining these two estimates. 
\end{Proof}

For later use, we will state the following lemma in a slightly 
more general form than 
what we will use in the proof of Theorem~\ref{th:main2}. 

\begin{Lem} \label{lem:time-3conti}
Let $( \mu_s )_{s \in [ 0 , \infty )}$ and $( \nu_s )_{s \in [ 0 , \infty )}$ 
be families of probability measures on $M$ which is continuous 
in $s$ with respect to the topology of weak convergence. 
For $t > 0$, 
let $\hat{d} : [ 0 , t ] \times M \times M \to [ 0 , \infty )$ 
be a continuous function such that 
$\hat{d} ( s , \cdot , \cdot )$ is a distance function on $M$ 
for each $s \in [ 0 , t ]$. 
Then
$s \mapsto \mathcal{T}_{\ph_{t-s} (\hat{d}(s, \cdot, \cdot))} ( \mu_s , \nu_s )$ 
is continuous on $[ 0 , t )$ and lower semi-continuous at $t$. 
\end{Lem}

\begin{Proof}
For simplicity of notations, 
we denote $\hat{d} ( s , x , y )$ by $d_s (x,y)$. 
Let $s_0 \in [ 0 , t )$ and 
take a decreasing sequence $( s_n )_{n \in \N}$ 
with $\lim_{n \to \infty} s_n = s_0$. 
Let $\ep > 0$. 
Since $( \mu_{s_n} )_{n \in \N}$ and $( \nu_{s_n} )_{n \in \N}$ 
are tight in $\mathcal{P} (M)$,  
there exist a compact set $K \subset M$ such that 
$\pi ( ( K \times K )^c ) < \ep$ 
for any $\pi \in \bigcup_{n \in \N} \Pi ( \mu_{s_n} , \nu_{s_n} )$. 
Since $\ph_{t - s_n} (a)$ is nonincreasing in $n$, 
the Dini theorem yields that 
$\ph_{t - s_n} ( d_{s_0} )$ converges to 
$\ph_{t - s_0} ( d_{s_0} )$ as $n \to \infty$ 
uniformly on $K \times K$. 
By Corollary~\ref{cor:distance}, 
\begin{equation*}
\abs{ 
  \ph_{t - s_n} ( d_{s_n} (x,y) ) 
  - 
  \ph_{t - s_n} ( d_{s_0} (x,y) )
}
\le 
\ph_{t - s_n} ( | d_{s_n} (x,y) - d_{s_0} (x,y) | ). 
\end{equation*}
By the assumption on $d_s$, we have 
\begin{equation*}
\lim_{n \to \infty} 
\sup_{x,y \in K} | d_{s_n} (x,y) - d_{s_0} (x,y) | 
= 0. 
\end{equation*} 
By combining these estimates, we obtain 
\begin{align} \nonumber
\limsup_{n \to \infty} 
& 
\abs{
  \mathcal{T}_{\ph_{t - s_n} ( d_{s_n} )} ( \mu_{s_n} , \nu_{s_n} )
  - 
  \mathcal{T}_{\ph_{t - s_0} ( d_{s_0} )} ( \mu_{s_n} , \nu_{s_n} )
} 
\\ \nonumber
& \le 
\limsup_{n \to \infty} 
\bigg( 
\sup_{x,y \in  K} 
\abs{ 
  \ph_{t - s_n} ( d_{s_0} ( x, y ) ) 
  - 
  \ph_{t - s_0} ( d_{s_0} ( x, y ) ) 
}
\\ & \hspace{9em} \nonumber 
+ 
\sup_{x,y \in K} 
\ph_{t - s_1} 
\abra{ | d_{s_n} ( x, y ) - d_{s_0} ( x, y ) | }
\bigg)
+ \ep
\\ \label{eq:3conti1}
& = \ep .
\end{align}
By virtue of Corollary~\ref{cor:distance} and 
\cite[Theorem~7.12]{book_Vil1}, 
the weak convergences 
$\mu_{s_n} \to \mu_{s_0}$ 
and 
$\nu_{s_n} \to \nu_{s_0}$ 
imply that 
$\mathcal{T}_{\ph_{s_0} ( d_{s_0} )} ( \mu_{s_n} , \nu_{s_n} )$ converges 
to $\mathcal{T}_{\ph_{s_0} ( d_{s_0} )} ( \mu_{s_0} , \nu_{s_0} )$. 
By combining this fact with \eqref{eq:3conti1}, 
we obtain 
\begin{equation*}
\lim_{n \to \infty} 
\mathcal{T}_{\ph_{t - s_n} ( d_{s_n} ) } ( \mu_{s_n} , \nu_{s_n} )
= 
\mathcal{T}_{\ph_{t - s_0} ( d_{s_0} ) } ( \mu_{s_0} , \nu_{s_0} ). 
\end{equation*}
It proves that 
$\mathcal{T}_{\ph_{t - s} ( d_{s} ) } ( \mu_{s} , \nu_{s} )$ 
is right-continuous at $s_0$. 
In a similar way, we can show the left-continuity of 
$\mathcal{T}_{\ph_{t - s} ( d_{s} ) } ( \mu_{s} , \nu_{s} )$ 
at $s_0$. 
Finally we will show the lower semi-continuity at $t$. 
Since $\ph_r (a)$ is nonincreasing in $r$, 
for $t' > t$, we have 
\begin{equation*}
\liminf_{s \uparrow t} \mathcal{T}_{\ph_{t-s} (d_s)} ( \mu_s , \nu_s ) 
\ge 
\lim_{s \uparrow t} \mathcal{T}_{\ph_{t' - s} (d_s)} ( \mu_s , \nu_s )
= 
\mathcal{T}_{\ph_{t'-t} (d_t)} ( \mu_t , \nu_t ). 
\end{equation*}
Hence the conclusion follows from 
Lemma~\ref{lem:time-Wconti} \ref{item:W-linf}
by letting $t' \downarrow t$. 
\end{Proof}

\begin{tProof}{\protect{Theorem~\ref{th:main2}}}
Recall that, for $t' > 0$, $s' \ge 0$ and $x_1 , x_2 \in M$,
Theorem~\ref{th:main} yields 
\begin{equation} \label{eq:pre-W0}
\mathcal{T}_{\ph_{t'} (d)} ( \P_{x_1} \circ X(s')^{-1} , \P_{x_2} \circ X (s')^{-1} ) 
\le 
\ph_{t' + s'} ( d ( x_1 , x_2 ) ). 
\end{equation}
Let $0 \le s_1 \le s_2 < t$. 
For each $y_1 , y_2 \in M$, 
take 
$
\pi_{s_2 - s_1}^{y_1 y_2}  
\in 
\Pi ( \P_{y_1} \circ X (s_2 - s_1)^{-1} , \P_{y_2} \circ X (s_2 - s_1 )^{-1} )
$ 
so that 
\begin{equation*}
\mathcal{T}_{\ph_{t - s_2} (d)} 
( 
  \P_{y_1} \circ X (s_2 - s_1 )^{-1} , 
  \P_{y_2} \circ X ( s_2 - s_1 )^{-1} 
)
= 
\int_{M \times M} \ph_{t - s_2} (d) 
\; d \pi_{s_2 - s_1}^{y_1 y_2}. 
\end{equation*}
We can choose $\pi_{s_2 - s_1}^{y_1 y_2}$ so that 
$( y_1 , y_2 ) \mapsto \pi_{s_2 - s_1}^{y_1 y_2}$ is measurable 
(see \cite[Corollary~5.22]{book_Vil2}, for instance).
Let us take a minimizer $\pi \in \Pi ( \mu^{(1)}_{s_1} , \mu^{(2)}_{s_1} )$ 
of $\mathcal{T}_{\ph_{t - s_1} (d)} ( \mu^{(1)}_{s_1} , \mu^{(2)}_{s_1} )$ 
and define $\pi^* \in \Pi ( \mu^{(1)}_{s_2} , \mu^{(2)}_{s_2} )$ 
by 
\[
\pi^* (A) := \int_{M \times M} \pi_{s_2 - s_1}^{y_1 y_2} (A) \pi ( dy_1 dy_2 ).  
\]
Then, 
by applying \eqref{eq:pre-W0} 
with $t' = t - s_2$ and $s' = s_2 - s_1$, 
we obtain 
\begin{align*}
\mathcal{T}_{\ph_{t-s_2} (d)} ( \mu^{(1)}_{s_2} , \mu^{(2)}_{s_2} ) 
& \le 
\int_{M \times M} \ph_{t - s_2} (d) \, d \pi^* 
\\
& = 
\int_{M \times M} 
\mathcal{T}_{\ph_{t-s_2} (d)} ( \P_{y_1} \circ X (s_2 - s_1 )^{-1} , \P_{y_2} \circ X (s_2 - s_1 )^{-1} ) 
\pi (dy_1 dy_2 )
\\
& \le 
\int_{M \times M} \ph_{t - s_1} ( d ( y_1 , y_2 ) ) \pi (dy_1 dy_2 )
\\
& = 
\mathcal{T}_{\ph_{t-s_1} (d)} ( \mu^{(1)}_{s_1} , \mu^{(2)}_{s_1} ). 
\end{align*}
Thus the assertion holds when $t > s_2$. When $t = s_2$, 
the assertion follows by taking $s_2 \uparrow t$ 
together with Lemma~\ref{lem:time-3conti} 
with $\hat{d} (t,x,y) : = d (x,y)$. 
\end{tProof}

In Theorem~\ref{th:main2}, the cost function $\ph_{t-s} (d)$ depends 
on time parameter $t$. Thus it seems to be natural to 
consider the limit $t \to \infty$, 
under a suitable scaling if necessary. 
We can realize it  
by combining Corollary~\ref{cor:t-large} with Theorem~\ref{th:main2} 
with the aid of 
Lemma~\ref{lem:time-Wconti}. 
Then we obtain the following monotonicity of transportation costs: 

\begin{Cor} \label{cor:mono}
Let us define 
$\Theta = \Theta_{K,N} : [ 0 , \bar{R} ) \to [ 0 , \infty )$ and 
$\kappa = \kappa (K,N) \in \R$
by 
\begin{align*}
\Theta_{K,N} (a) 
& : = 
\begin{cases} 
a 
& 
(\mbox{$K = 0$}), 
\\
a 
& 
(\mbox{$N = \infty$ and $K > 0$}),
\\
\displaystyle
\chi \abra{ \frac{ a \sqrt{-K}}{2} } 
& 
(\mbox{$N = \infty$ and $K < 0$}), 
\\
\displaystyle
\sin \abra{ \frac12 \sqrt{\frac{K}{N-1}} a}
& 
(\mbox{$N < \infty$ and $K > 0$}), 
\\
\displaystyle
\int_0^\infty 
\chi \abra{ 
  \sqrt{\frac{-Ku}{2(N-1)}} 
  s_{K/(N-1)} \abra{\frac{a}{2}}
}
u^{(N-3)/2} \e^{-u} du 
& 
(\mbox{$N < \infty$ and $K < 0$}),  
\end{cases}
\\
\kappa (K,N) 
& := 
\begin{cases}
K \vee 0 
& (N = \infty), 
\\
\displaystyle 
\frac{NK}{N-1} \vee 0 
& ( N < \infty ). 
\end{cases}
\end{align*} 
For $i=1, 2$ and $\mu^{(i)} \in \mathcal{P} (M)$, 
let $\mu^{(i)}_t$ be the distribution of $X(t)$ 
with the initial distribution $\mu^{(i)}$. 
Then 
$\e^{\kappa s} \mathcal{T}_{\Theta (d)} ( \mu^{(1)}_s , \mu^{(2)}_s )$ 
is nonincreasing in $s$. 
\end{Cor}

When $K = 0$ or $N = \infty$ and $K > 0$, 
what the last corollary states is nothing but 
the $L^1$-Wasserstein contraction. 
When $N < \infty$ and $K > 0$, 
what we obtained is essentially well-known 
(see \cite{Wang_equivCD} for the statement 
formulated in terms of optimal transport theory). 
Thus the most interesting assertion is in the case $K < 0$. 
In the usual $L^p$-Wasserstein contraction in \eqref{eq:W-contr}, 
The upper bound grows exponentially fast as time increases 
when $K < 0$. 
The last corollary says that a nonincreasing property still holds 
even when $K < 0$ by choosing a cost function appropriately. 

\section{Gradient estimates}
\label{sec:gradient}

For a bounded and measurable function $f : M \to \R$, 
we define the action of the diffusion semigroup $P_t f$ 
by $P_t f (x) : = \E_x [ f ( X ( t ) ) ]$. 
We denote the dual action of $P_t$
to $\mathcal{P} (M)$ by $P_t^*$. 
That is, 
\[
P_t^* \mu (A) = \int_M \P_x [ X (t) \in A ] \, \mu (dx).
\]
Since $\ph_t (d(x,y))$ is a distance function 
by Corollary~\ref{cor:distance}, the Kantorovich-Rubinstein 
duality easily implies the following 
(cf.~\cite{K9,Sturm_SemigrHarm}): 

\begin{Thm} \label{th:gradient}
Given $t, s \ge 0$, the following are equivalent: 
\begin{enumerate}
\item \label{item:T}
For $\mu_1 , \mu_2 \in \mathcal{P}(M)$, 
\[
\mathcal{T}_{\ph^{K,N}_t (d)} ( P_s^* \mu_1 , P_s^* \mu_2 )
\le
\mathcal{T}_{\ph^{K,N}_{t+s} (d)} ( \mu_1 , \mu_2 ).
\]
\item \label{item:G}
For any $\ph^{K,N}_t (d)$-Lipschitz function $f$ on $M$, 
\[
\sup_{x\neq y} 
\frac{| P_s f (x) - P_s f (y) |}{\ph^{K,N}_{t+s} ( d (x,y) )} 
\le 
\sup_{x \neq y} \frac{| f(x) - f(y) |}{\ph^{K,N}_t (d(x,y))}. 
\]
\end{enumerate}
\end{Thm}
The condition \ref{item:T} in the last theorem 
comes from the consequence of Theorem~\ref{th:main2}. 
Note that, in the condition \ref{item:G}, 
those supremums may be attained 
at $(x,y) \in M \times M$ with $d (x,y) > 0$ 
since $\ph_t (d)$ is not a geodesic distance. 
As an easy consequence of Theorem~\ref{th:gradient}, 
we obtain the following gradient estimate. 

\begin{Cor} \label{cor:gradient}
Under Assumption~\ref{ass:CD}, we have 
\begin{equation*}
\| \nabla P_t f \|_\infty 
\le 
\ph_t' (0) \osc (f) 
\end{equation*}
for any bounded measurable function $f$ on $M$. 
\end{Cor}
Recall that an expression of $\ph_t'(0)$ is given in Proposition~\ref{prop:concave}. 
Note that a gradient estimate like in Corollary~\ref{cor:gradient} 
also follows from the reverse Poincar\'e inequality 
(see \cite{Bak97,Led_geom-Markov} for instance; also see \cite{AGS3}). 
When $K \ge 0$, Corollary~\ref{cor:gradient} and 
Proposition~\ref{prop:concave} \ref{item:diff_at_zero}
easily imply the Liouville property, that is, 
there are no nonconstant bounded $\mathcal{L}$-harmonic functions, 
by taking $f$ as a bounded harmonic function (so that $P_t f = f$) 
and $t \to \infty$. 
\begin{Proof}
Theorem~\ref{th:main2} tells us that the condition \ref{item:T}
holds with $t = 0$ under Assumption~\ref{ass:CD}. 
Then the definition of $\ph_0$ yields 
\[
\sup_{x \neq y} \frac{| f(x) - f(y) |}{\ph^{K,N}_0 (d(x,y))} = \osc (f). 
\]
Recall that the differentiability of $\ph_s^{K,N}$ at $0$
is given in 
Proposition~\ref{prop:concave} \ref{item:diff_at_zero}. 
Thus Theorem~\ref{th:gradient} implies that 
\[
\frac{1}{(\ph^{K,N}_s)' (0)}\| \nabla P_s f \|_\infty 
\le 
\sup_{x \neq y} 
\frac{ | P_s f (x) - P_s f (y) |}{\ph_s (d(x,y))}
\le \osc (f)
\]
and hence the conclusion holds. 
\end{Proof}

\section{Stability under the Gromov-Hausdorff convergence}
\label{sec:stable}

In this section we consider a sequence of 
Riemannian manifolds $( M_n , g_n )$ ($n \in \N$). 
By technical reasons, 
we restrict ourselves into the case that 
each $M_n$ is compact. 
Let $f_n$ be a positive $C^1$-function on $M_n$ 
and $Z_n : = \nabla f_n$ for $n \in \N$. 
Suppose that, given $N < \infty$ and $K$, 
$( M_n , g_n )$ and $Z_n$ satisfies Assumption~\ref{ass:CD} 
for all $n \in \N$ where the parameters $N$ and $K$ are independent of $n$. 
Let $\mathrm{vol}_n$ be the Riemannian volume measure 
on $( M_n , g_n )$ and set $\nu_n = \e^{f_n} \mathrm{vol}_n$. 
Under Assumption~\ref{ass:CD}, the metric measure space $( M_n , d_n , \nu_n )$ 
satisfies the curvature-dimension condition $\mathsf{CD} ( K , \infty )$ 
(see \cite{Lott-Vill_AnnMath09,Sturm_Ric,Sturm_Ric2}). 
Thus the gradient flow of the relative entropy functional 
$\mathrm{Ent}_{\nu_n}$ on $L^2$-Wasserstein space 
over $( M_n , d_n , \nu_n )$ 
is identified with the gradient flow of 
the Dirichlet energy functional on $L^2 ( M_n , \nu_n )$ 
(see \cite{AGS2,Erb_gradF-Riem,GKO}). 

\begin{Defn} \label{def:eps-isom}
\begin{enumerate}
\item
Let $( M_1 , d_1 )$ and $( M_2 , d_2 )$ be 
metric spaces. 
For $\ep > 0$, we call $f : M_1 \to M_2$ an $\ep$-isometry 
if the following hold: 
\begin{align*}
\sup_{x , y \in M_1} | d_1 ( x , y ) - d_2 ( f (x) , f (y) ) | 
& \le \ep , 
& \sup_{y \in M_2 } d_2 ( y , f ( M_1 ) ) 
& \le \ep .
\end{align*}
\item
Let $( ( M_n , d_n , \nu_n ) )_{n \in \N}$ and 
$( M , d , \nu )$ be metric measure spaces. 
We say $( M_n , d_n , \nu_n )$ converges to 
$( M , d , \nu )$ as $n \to \infty$ in 
the measured Gromov-Hausdorff sense 
if there exist $\ep_n > 0$ ($n \in \N$) 
with $\lim_{n \to \infty} \ep_n = 0$ and 
$\ep_n$-isometry $f_n : M_n \to M$ so that 
$f_n^\# \nu_n$ converges to $f^\# \nu$ 
in the vague topology. 
\end{enumerate}
\end{Defn}

In the sequel, we assume that 
$( M_n , d_n , \nu_n )$ converges to 
a compact metric measure space $( M , d, \nu )$ 
in the measured Gromov-Hausdorff sense 
via $\ep_n$-isometries $f_n : M_n \to M$. 
Note that, in this framework, 
the convergence with respect to 
the measured Gromov-Hausdorff distance 
is equivalent to the convergence with respect to 
the distance $\mathbf{D}$ introduced in 
\cite{Sturm_Ric} (see \cite[Subsection~3.4]{Sturm_Ric}). 
Under the assumption, 
$( M, d , \nu )$ satisfies $\mathsf{CD} (K , \infty )$ 
again (see \cite{AGS3,Lott-Vill_AnnMath09,Sturm_Ric}). 
Thus, for $\mu_0 \in \mathcal{P} (M)$ 
with $\mathrm{Ent}_\nu ( \mu_0 ) < \infty$, 
there exists a unique gradient curve $\mu_t$ 
of $\mathrm{Ent}_\nu$ on $\mathcal{P} ( M )$ 
starting from $\mu_0$ 
(see \cite{AGS2,AGS3,Gigli_Heatmm}). 

The following theorem asserts that 
these gradient curves enjoy the same monotonicity as 
shown in Theorem~\ref{th:main2}: 

\begin{Thm} \label{th:stability}
For $i=1, 2$, let $\mu^{(i)}_0 \in \mathcal{P} (M)$ 
with $\mathrm{Ent}_\nu ( \mu^{(i)}_0 ) < \infty$ 
and $\mu^{(i)}_t$ a gradient curve of $\mathrm{Ent}_{\nu}$ 
with initial distribution $\mu^{(i)}_0$. 
Then, for any $t \in [ 0 , \infty )$, 
$\mathcal{T}_{\ph_{t - s} ( d )} ( \mu^{(1)}_s , \mu^{(2)}_s )$ is 
a nonincreasing function of $s \in [ 0 , t ]$.  
\end{Thm}

\begin{Proof}
By virtue of Lemma~\ref{lem:time-3conti}, 
it suffices to show the assertion for $s \in (0,t)$. 
For $i = 1, 2$, 
there exists $\mu^{(i,n)}_0 \in \mathcal{P} ( M_n )$ 
for $n \in \N$ 
such that $\mathrm{Ent}_{\nu_n} ( \mu^{(i,n)}_0 ) < \infty$ and 
that $f_n^\# \mu^{(i,n)}_0$ converges to $\mu^{(i)}_0$ 
by following an argument 
in the proof of \cite[Theorem~4.15]{Lott-Vill_AnnMath09}. 
Let $\mu^{(i,n)}_t$ be the gradient curve of $\mathrm{Ent}_{\nu_n}$ 
on $\mathcal{P} ( M_n )$ 
with the initial distribution $\mu^{(i,n)}_0$. 
Then, by virtue of \cite[Theorem~21]{Gigli_Heatmm}, 
$f_n^\# \mu^{(i,n)}_t$ converges to $\mu^{(i)}_t$ for $i = 1, 2$ and $t > 0$. 

We claim that for each $s \in ( 0, t )$, 
\begin{equation} \label{eq:mGH0}
\lim_{n \to \infty} 
\abra{ 
  \mathcal{T}_{\ph_{t-s} (d_n)} ( \mu^{(1,n)}_s , \mu^{(2,n)}_s )
  - 
  \mathcal{T}_{\ph_{t-s} (d)} ( f_n^\# \mu^{(1,n)}_s , f_n^\# \mu^{(2,n)}_s )
} 
= 0. 
\end{equation}
Take $\pi^{(n)} \in \Pi ( \mu^{(1,n)}_s , \mu^{(2,n)}_s )$ and 
set $\tilde{\pi}^{(n)} := \abra{ f_n \times f_n }^\# \pi^{(n)}$. 
Then we can easily see that 
$\tilde{\pi}^{(n)} \in \Pi ( f_n^\# \mu^{(1,n)}_s , f_n^\# \mu^{(2,n)}_s )$ 
holds. 
Since $f_n$ is $\ep_n$-isometry and 
$\ph_{t-s} ( \cdot )$ is nondecreasing, 
\begin{align} \nonumber
\mathcal{T}_{\ph_{t-s} (d)} ( f_n^\# \mu^{(1,n)}_s , f_n^\# \mu^{(2,n)}_s )
& \le 
\int_{M \times M} 
 \ph_{t-s} ( d (x,y) ) 
\tilde{\pi}^{(n)} ( dx dy ) 
\\ \nonumber
& = 
\int_{M \times M} 
 \ph_{t-s} ( d ( f_n (x) , f_n (y) ) ) 
\pi^{(n)} ( dx dy )
\\ \label{eq:mGH1}
& = 
\int_{M \times M} 
 \ph_{t-s} ( d_n (x,y) + \ep_n ) 
\pi^{(n)} ( dx dy )
.  
\end{align}
Since our choice of 
$\pi^{(n)} \in \Pi ( \mu^{(1,n)}_s , \mu^{(2,n)}_s )$
can be arbitrary, 
\eqref{eq:mGH1} and Corollary~\ref{cor:distance} yield 
\begin{equation} \label{eq:mGH2}
\mathcal{T}_{\ph_{t-s} (d)} ( f_n^\# \mu^{(1,n)}_s , f_n^\# \mu^{(2,n)}_s ) 
\le 
\mathcal{T}_{\ph_{t-s} (d_n)} ( \mu^{(1,n)}_s , \mu^{(2,n)}_s ) 
+ 
\ph_{t-s} ( \ep_n ) 
. 
\end{equation}
To complete the proof of the claim, 
let us take an approximate inverse $g_n$ for each $n \in \N$, 
that is, $g_n : M \to M_n$ satisfies 
\begin{align*}
\lim_{n \to \infty} 
& 
\sup_{ x \in M_n } 
d ( x , g_n ( f_n (x) ) ) 
= 0 ,
& 
\lim_{n \to \infty} 
& 
\sup_{ x \in M } 
d_n ( x, f_n ( g_n (x) ) )
= 0 . 
\end{align*}
We may assume $g_n$ is $\ep_n'$-isometry 
for some $\ep_n'$ with $\lim_{n \to \infty} \ep_n' = 0$ 
without loss of generality. 
By a similar argument as what we used to obtain \eqref{eq:mGH2}, 
\begin{multline} \label{eq:mGH3}
\mathcal{T}_{\ph_{t-s} (d_n)} 
\abra{ 
  ( g_n \circ f_n )^\# \mu^{(1,n)}_s 
  , 
  ( g_n \circ f_n )^\# \mu^{(2,n)}_s 
} 
\\
\le 
\mathcal{T}_{\ph_{t-s} (d)} 
\abra{
  f_n^\# \mu^{(1,n)}_s 
  , 
  f_n^\# \mu^{(2,n)}_s 
} 
+ 
\ph_{t-s} ( \ep_n' ) 
. 
\end{multline}
Since 
$
( \mathrm{id} \times (g_n \circ f_n ) )^\# \mu^{(i,n)}_s
\in 
\Pi ( \mu^{(i,n)}_s , ( g_n \circ f_n )^\# \mu^{(i,n)}_s )
$, 
\begin{align*}
\mathcal{T}_{\ph_{t-s} (d_n)} 
\abra{
  \mu^{(i,n)}_s 
  , 
  ( g_n \circ f_n )^\# \mu^{(i,n)}_s 
} 
& \le 
\int_{M_n \times M_n} 
d_n ( x , g_n ( f_n (x) ) ) \mu^{(i,n)} (dx)
\\
& \le 
\sup_{x \in M_n} d_n ( x, g_n ( f_n (x) ) ) 
\end{align*}
for $i = 1,2$. 
By combining this estimate with \eqref{eq:mGH3}, 
we obtain 
\begin{multline} \label{eq:mGH4}
\mathcal{T}_{\ph_{t-s} (d_n)} 
\abra{ 
  \mu^{(1,n)}_s 
  , 
  \mu^{(2,n)}_s 
} 
\\
\le 
\mathcal{T}_{\ph_{t-s} (d)} 
\abra{
  f_n^\# \mu^{(1,n)}_s 
  , 
  f_n^\# \mu^{(2,n)}_s 
} 
+ 
2 \sup_{x \in M_n} d_n ( x, g_n ( f_n (x) ) ) 
+ 
\ph_{t-s} ( \ep_n' ) 
. 
\end{multline}
Hence \eqref{eq:mGH2} and \eqref{eq:mGH4} imply the claim 
since $\ph_{t-s} (\cdot)$ is continuous. 

By Corollary~\ref{cor:distance} and \cite[Theorem~7.12]{book_Vil1}, 
$\mathcal{T}_{\ph_{t-s} (d)} ( f_n^\# \mu^{(i,n)}_s , \mu^{(i)}_s )$ 
converges to 0 as $n \to \infty$ for $i = 1,2$. 
Hence \eqref{eq:mGH0} yields 
$
\lim_{n \to \infty} 
\mathcal{T}_{\ph_{t-s} ( d_n ) } ( \mu^{(1,n)}_s , \mu^{(2,n)}_s ) 
= 
\mathcal{T}_{\ph_{t-s} (d)} ( \mu^{(1)}_s , \mu^{(2)}_s )
$. 
Since $\mathcal{T}_{\ph_{t-s} (d_n)} ( \mu^{(1,n)}_s , \mu^{(2,n)}_s )$ 
is nonincreasing in $s$ by Theorem~\ref{th:main2}, 
the conclusion holds.  
\end{Proof}

\section{Time-dependent metrics}

Let $( g(t) )_{t \in [ T_1 , T_2 ]}$ be a family of 
smooth complete Riemannian metrics on $M$ depending 
smoothly in $t$. 
Let $Z (t)$ be a time-dependent vector field on $M$ 
depending continuously in $t$ and consider 
the time-inhomogeneous diffusion process 
$( ( X (t) )_{t \in [ T_1 , T_2 ]} , ( \P_x )_{x \in M} )$ 
generated by $\mathcal{L}_t := \Dl_{g(t)} + Z (t)$. 
The following assumption corresponds to 
Assumption~\ref{ass:CD} with $N = \infty$: 
\begin{Ass} \label{ass:tCD}
Given $K \in \R$, 
the following holds for each $t$: 
\[
(\nabla Z(t))^\flat + \frac12 \partial_t g(t) 
\le 
\Ric_{g(t)} - K g(t). 
\]
\end{Ass}
An important example of the time-dependent metrics $g(t)$ 
satisfying Assumption~\ref{ass:tCD} is the backward 
Ricci flow, that is, 
\[
\frac12 \partial_t g(t) 
= \Ric_{g(t)}. 
\]
Under Assumption~\ref{ass:tCD}, 
the coupling by reflection of $X(t)$ is already 
studied in \cite{K10} via the approximation 
by geodesic random walks 
(The notation in \cite{K10} looks slightly different 
since we considered the diffusion process 
generated by $\Dl_{g(t)} / 2  + Z (t)$ there). 
By modifying arguments in previous sections, 
we can obtain the results corresponding to 
Theorem~\ref{th:main}, Theorem~\ref{th:main2}, 
Corollary~\ref{cor:TV-comparison} and Corollary~\ref{cor:mono} 
with $N = \infty$ by replacing 
$d$ which measures the distribution at $s$ with $d_{g ( T_1 +s )}$. 
For example, the conclusion of the statement corresponding to 
Theorem~\ref{th:main2} is as follows: Let $\mu^{(i)}_s$ be 
the distribution of $X (t)$ at $t = s + T_1$ with initial 
distribution $\mu^{(i)}_{T_1}$ for $i = 1,2$. 
Then, for $t \ge s_2 > s_1 \ge 0$, 
\begin{equation} \label{eq:t-mono}
\mathcal{T}_{\ph_{t-s_2} ( d_{g ( s_2 + T_1 )} )} 
( \mu^{(1)}_{s_2 + T_1} , \mu^{(2)}_{s_2 + T_1} ) 
\le
\mathcal{T}_{\ph_{t-s_1} ( d_{g ( s_1 + T_1 )} )} 
( \mu^{(1)}_{s_1 + T_1} , \mu^{(2)}_{s_1 + T_1} ). 
\end{equation} 

For reader's convenience, let us explain briefly 
why the time derivative with respect to the metric 
appears in Assumption~\ref{ass:tCD}. 
When we follow the argument in the time-independent metric case 
in Proposition~\ref{prop:DI}, 
we consider $d_{g(t^\a_n)} ( \mathbf{X}^\a (t^\a_n) )$ instead of 
$r^\a (n) = d (\mathbf{X}^\a ( t^\a_n ) )$. 
Then, in the Taylor expansion in the proof of Lemma~\ref{lem:c-2var}, 
there appears the time derivative of $d_{g(t)}$ as an additional term. 
It can be described in terms of the time derivative of $g(t)$. 
Then our condition in Assumption~\ref{ass:tCD} will be used to implement 
this additional term into the lower bound of Bakry-\'Emery Ricci tensor. 
For more details, see \cite{K10}. 

Note that, by \cite[Lemma~2.5]{K10}, 
$\hat{d} (t,x,y) := d_{g(t + T_1)} (x,y)$ satisfies 
the assumption of Lemma~\ref{lem:time-3conti}. 
This fact will be used to complete the proof of \eqref{eq:t-mono} 
when $t = s_2$. 
\bigskip

\noindent 
\textit{Acknowledgment.} 
The first named author is grateful for the support 
by the Grant-in-Aid for Young Scientists (B) 22740083. 

\bibliographystyle{amsplain}
\providecommand{\bysame}{\leavevmode\hbox to3em{\hrulefill}\thinspace}
\providecommand{\MR}{\relax\ifhmode\unskip\space\fi MR }
% \MRhref is called by the amsart/book/proc definition of \MR.
\providecommand{\MRhref}[2]{%
  \href{http://www.ams.org/mathscinet-getitem?mr=#1}{#2}
}
\providecommand{\href}[2]{#2}

\end{document}